\documentclass[11pt,twoside]{article}

\usepackage{latexsym,amsmath,graphics}
\usepackage[psamsfonts]{amssymb}
\usepackage[bookmarks,bookmarksnumbered,pdfstartview=FitBH]{hyperref}

\newtheorem{theoremintro}{\textrm{Theorem}}[section]
\newtheorem{corollaryintro}[theoremintro]{\textrm{Corollary}}
\newtheorem{theorem}{Theorem}[section]

\newtheorem{lemma}[theorem]{Lemma}
\newtheorem{proposition}[theorem]{Proposition}
\newtheorem{corollary}[theorem]{Corollary}
\newtheorem{definition}[theorem]{Definition}
\newtheorem{example}[theorem]{Example}
\newtheorem{remark}[theorem]{Remark}

\newcommand{\scal}[2]{{\langle #1 \vert #2 \rangle}}

\def\RR{\mathcal{R}}

\def\FF{\mathbf{F}}
\def\Nmath{\mathbb{N}}
\def\Zmath{\mathbb{Z}}
\def\Rmath{\mathbb{R}}
\def\Cmath{\mathbb{C}}

\def\proof{\textsc{Proof:\ }}
\def\endofproof{\nobreak\hfill $\blacksquare$ \goodbreak}  

\def\GG{\mathcal{G}}  
\def\EE{\mathtt{E}}  
\def\VV{\mathtt{V}}  
\def\GGamma{\Gamma}
\def\GGG{G}  
\def\HH{H} 

\def\origin{\rho} 
\def\Xbar{\overline{X}}
\def\mubar{\overline{\mu}}
\def\RRbar{\overline{\mathcal{R}}}
\def\full{\mathrm{fu}}
\def\cluster{\mathrm{cl}}
\def\Lfull{\mathcal{L}^{\full}}
\def\RRfull{\mathcal{R}^{\full}}
\def\Lcluster{\mathcal{L}^{\cluster}}
\def\RRcluster{\mathcal{R}^{\cluster}}
\def\HD{\mathbf{HD}}
\def\OOO{\mathcal{O}_{\mathbf{HD}}}
\def\HHH{\bar{H}_{\!(2)}^{1}}
\def\dzero{d}
\def\dun{D}
\def\dim{\mathrm{dim}}

\def\Omegaunique{U} 
\def\Xunic{U^{\pi}}
\def\hh{h}
\def\KK{K_{\origin}} 
\def\kk{k}
\def\Cbar{\overline{C}}
\def\Haar{m} 
\def\WW{W}

\def\TT{T^{\bullet}}

\def\RC{\mathtt{RC}}
\def\FRC{\mathtt{FRC}}
\def\WRC{\mathtt{WRC}}
\def\val{\mathrm{deg}}

\begin{document}

\thispagestyle{empty}

\date{}
\title{Invariant Percolation and Harmonic Dirichlet Functions}
\author{Damien Gaboriau\thanks{C.N.R.S.}}

\maketitle

\begin{abstract} The main goal of this paper is to  answer question~1.10 and 
settle conjecture~1.11 of Benjamini-Lyons-Schramm \cite{BLS99} relating harmonic 
Dirichlet functions on a graph to those on the infinite clusters in the 
uniqueness phase of Bernoulli percolation. We extend the result to more general 
invariant percolations, including the Random-Cluster model. We prove the existence 
of the nonuniqueness phase for the Bernoulli percolation (and make some progress 
for Random-Cluster model) on unimodular transitive locally finite graphs admitting 
nonconstant harmonic Dirichlet functions. This is done by using the device of 
$\ell^2$ Betti numbers.
\end{abstract}


\addtocounter{section}{-1}

\section{Introduction}

Traditionally, percolation on graphs lives on 
$\Zmath^{d}$ or lattices in $\Rmath^{d}$.
Following earlier work of G.~Grimmett and C.~Newman \cite{GN90} on the direct
product of a regular tree and $\Zmath$, a general study of
invariant percolation was initiated in I.~Benjamini and O.~Schramm~\cite{BS96}
and further developed by several authors.

\medskip
Let $\GG=(\VV,\EE)$ be a (non-oriented) countable infinite locally finite graph. 
A \textbf{bond percolation} 
on $\GG$ is simply a probability measure $\mathbf{P}$ on $\Omega=\{0,1\}^{\EE}$, the subsets of its edge set $\EE$.
It is an \textbf{invariant percolation} when this 
measure is invariant under a certain group $\HH$ of automorphisms of $\GG$.

An element $\omega$ in $\Omega$ defines the \textbf{graph} whose vertices are 
$\VV$ and whose edges are the \textbf{retained} (or \textbf{open}) 
edges, {i.e.\ } those $e\in \EE$ with value $\omega(e)=1$.
It is the subgraph of $\GG$ where edges with value $0$ are 
\textbf{removed} (or \textbf{closed}).
One is interested in the shape of the ``typical'' random 
subgraph $\omega$\,\footnote{In more probabilistic terms, $\omega$ is a random 
variable with values in $\Omega$ and distribution $\mathbf{P}$.} and of 
 its \textbf{clusters}, {i.e.\ } its connected components.

\medskip
One of the most striking instances is \textbf{Bernoulli bond percolation}, and 
particularly on a \textbf{Cayley graph}\,\footnote{A Cayley graph will 
    always be assumed to be for a finitely generated group and with respect to a {\em finite}
    generating system.}
of a finitely generated group:
each edge of $\GG$ is removed with probability 
$1-p$ independently (where $p\in[0,1]$ is a parameter). The resulting 
probability measure $\mu_{p}$ on $\Omega$ is the 
product Bernoulli measure $(1-p,p)$ on $\{0,1\}$. It is 
invariant under every automorphism group of $\GG$.
How does the behavior evolve as $p$ varies?
For small $p$, the clusters are a.s.\ all finite, while for $p=1$ the 
measure concentrates on the infinite subgraph $\GG$ itself. 
Depending on the value of the parameter, $\mu_{p}$-almost 
every subgraph $\omega\in \Omega$ has
no infinite cluster, infinitely many infinite clusters (\textbf{nonuniqueness phase}) 
or only one infinite cluster (\textbf{uniqueness phase}). According 
to a somewhat surprising result of O.~H\"aggstr\"om and Y.~Peres 
\cite{HP99}, the 
phases are organized around two phase transitions for two critical values 
of $p$ depending on the graph $0<p_{c}(\GG)\leq p_{u}(\GG)\leq 1$, 
as summarized\,\footnote{For Cayley graphs, say. More generally this 
picture appears for Bernoulli percolation of unimodular quasi-transitive (see below for the definitions) graphs 
\cite{HP99}. R.~Schonmann has even removed the unimodularity assumption.} in the following picture: 
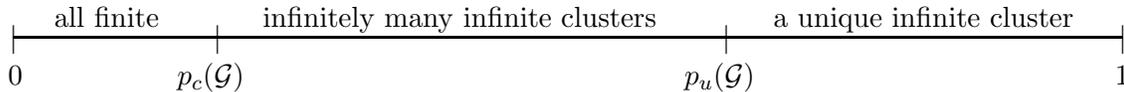
\begin{center}                                                               
\begin{tabular}{lcccccr}
    & \raisebox{-2ex}[0cm][0cm]{all finite} & & \hskip-5pt \raisebox{-2ex}[0cm][0cm]{infinitely many infinite 
    clusters} \hskip-5pt  & & \hskip-5pt  \raisebox{-2ex}[0cm][0cm]{a unique infinite cluster} \hskip-5pt &\\
    \multicolumn{7}{l}{\hrulefill}\hskip 2.5pt\\
    \hskip -1.5pt \raisebox{2.2ex}[0cm][0cm]{$\vert$} & & 
    \raisebox{2.2ex}[0cm][0cm]{$\vert$} &  &  \raisebox{2.2ex}[0cm][0cm]{$\vert$}
    & &  \raisebox{2.2ex}[0cm][0cm]{$\vert$}\phantom{\hskip 
    1pt}\\
    \hskip -2pt \raisebox{2ex}[0cm][0cm]{$0$} & & \hskip-5pt  \raisebox{2ex}[0cm][0cm]{$p_{c}(\GG)$} &  &  \hskip-5pt  
    \raisebox{2ex}[0cm][0cm]{$p_{u}(\GG)$} & & \raisebox{2ex}[0cm][0cm]{$1$}
\end{tabular}
\end{center}

The picture at the critical values themselves is far from complete
(to which interval belong the transitions? which inequalities are 
strict: $p_{c}\not= p_{u}\not= 1$ ?)
and seems to depend heavily (for Cayley graphs) on the algebraic properties of the group.
However, a certain amount of results has been obtained. For instance, 
in the Cayley graphs setting\,\footnote{{i.e.\ } for any Cayley graph $\GG$ 
of a given finitely generated group $\GGamma$}:

	\noindent $\bullet$ $p_{u}=p_{c}$ for amenable groups (Burton-Keane
	\cite{BK89})
    
	\noindent $\bullet$ $p_{c}<1$ for groups of polynomial or exponential 
    growth, except for groups with two ends \cite{Lyo95, LP05}

    \noindent $\bullet$  For any nonamenable group, there is 
    almost surely no infinite cluster at $p=p_{c}$ 
	\cite[Th.~1.3]{BLPS99a}\,\footnote{This is true more generally for nonamenable 
	unimodular transitive graphs \cite{BLPS99b}.}

    \noindent $\bullet$  $p_{u}<1$ for finitely presented groups with 
    one end (Babson-Benjamini \cite{BB99})
     and for (restricted) wreath products\,\footnote{
	The finitely generated group $\Lambda$ acts 
	transitively on the discrete set $\WW$ of indices and thus on the $\WW$-indexed direct product
	$\oplus_{\WW} K$.}
	$K\wr\Lambda:=\Lambda\ltimes \oplus_{\WW} K$
	with finite non-trivial $K$ (Lyons-Schramm \cite{LS99})

    \noindent $\bullet$  $p_{u}=1$ for groups with infinitely many 
    ends, thus the percolation 
    at $p=p_{u}$ belongs to the uniqueness phase\,\footnote{This is a very general 	result, see \cite{LP05}.}

    \noindent $\bullet$
    The percolation at the threshold $p=p_{u}$ belongs to the nonuniqueness 
    phase, and thus $p_{u}<1$, for infinite groups with Kazhdan's property (T)
    (Lyons-Schramm \cite{LS99})
    
    \noindent $\bullet$ in the nonuniqueness phase, infinite clusters have 
    uncountably many ends almost surely \cite{HP99}\,\footnote{This is 
    true more generally for quasi-transitive unimodular graphs \cite{HP99}.
	This has also been proved in the nonunimodular case by O.~H\"aggstr\"om, 
	Y.~Peres, R.~Schonmann \cite{HPS99}.}
\bigskip
\bigskip

For (much !)  more information and references, the reader is referred
to the excellent survey of R.~Lyons \cite{Lyo00}, book (in preparation)
by R.~Lyons and Y.~Peres \cite{LP05} and papers \cite{BLPS99a, BLPS99b, BLS99, BS96, HP99,LS99}.

\subsection{On Harmonic Dirichlet functions}

The space $\HD(\GG)$ of \textbf{Harmonic Dirichlet functions} on a
locally finite graph $\GG=(\VV,\EE)$ is the space of functions on the
vertex set $\VV$ whose value at each vertex equals the average of the
values at its neighbors $$\left(\displaystyle{\sum_{v'\sim v}}
1\right) f(v)=\displaystyle{\sum_{v'\sim v}} f(v')$$
and whose coboundary is $\ell^{2}$-bounded $$\Vert\dzero
f\Vert^{2}=\sum_{v\sim v'} \vert f(v)-f(v') \vert^{2}<\infty.$$
The constant functions on the vertex set $\VV$ always belong to
$\HD(\GG)$.  Denote by $\OOO$ the class of connected graphs for which these are
the only harmonic Dirichlet functions.  Belonging or not to $\OOO$
plays a role in electrical networks theory: when assigning resistance
1 ohm to each edge, the coboundary of a harmonic Dirichlet function
gives a finite energy current satisfying both Kirchhoff's laws.

As an
example, a Cayley graph of a group $\GGamma$ is in $\OOO$ if and
only if the first $\ell^{2}$ Betti number $\beta_{1}(\GGamma)$ of the
group vanishes (see Theorem~\ref{thm: beta-1 and HD of a Cayley graph}).
Thus, the Cayley graphs of the following groups (when finitely generated) 
all belong to $\OOO$: abelian groups,
amenable groups, groups with Kazhdan property (T), 
lattices in $\mathrm{SO}(n,1)$ ($n\geq 3$) or in 
$\mathrm{SU}(n,1)$.
 On the other hand, the class of groups whose Cayley graphs don't 
 belong to $\OOO$ contains the non-cyclic free groups, 
the fundamental groups of surfaces of genus $g\geq 2$, the 
free products of infinite groups, and the amalgamated free products
over an amenable group of groups in that class.
Look at the very informative paper by Bekka-Valette 
\cite{BV97} and F.~Martin's thesis \cite{Mar03} for further interpretations
in cohomological terms.

P.~Soardi \cite{Soa93} has proved that belonging to $\OOO$ is
invariant under a certain kind of ``perturbation'' of $\GG$, namely
quasi-isometry or rough isometry.  Bernoulli bond percolation clusters
may also be considered as perturbations of $\GG$. I.~Benjamini,
R.~Lyons and O.~Schramm addressed the analogous invariance 
problem\,\footnote{Observe that for $p_{u}<p < 1$, the infinite 
clusters are $\mu_p$ a.s.\ not quasi-isometric to $\GG$. They contain for instance 
arbitrarily long arcs without branch points (by deletion tolerance !).} in
\cite{BLS99}, by taking a stand only in one case:
\begin{quote}
    \textbf{Question \cite[Quest. 1.10]{BLS99}} Let $\GG$ be a Cayley graph, and suppose that $\GG\in \OOO$. 
    Let $\omega$ be Bernoulli percolation on $\GG$ in the uniqueness 
    phase. Does it follow that a.s.\ the infinite cluster of $\omega$ 
    is in $\OOO$?

\textbf{Conjecture \cite[Conj. 1.11]{BLS99}} Let $\GG$ be a Cayley graph, $\GG\not\in \OOO$. 
Then a.s.\ all infinite clusters of
$p$-Bernoulli percolation are not in $\OOO$.
\end{quote}
They proved this for $p$ sufficiently large:
\begin{quote}
    \textbf{Theorem}\textbf{\cite[Th. 1.12]{BLS99}} 
\emph{If a Cayley graph $\GG$ is not in $\OOO$, then there exists a $p_{0}<1$, such that every infinite cluster of $\mu_{p}$-a.e.\ subgraph is not in $\OOO$, for every $p\geq p_{0}$}.
\end{quote}
On the other hand, 
\begin{quote}
    \textbf{Theorem}\textbf{\cite[Th. 1.9]{BLS99}} \emph{If $\GG$ is a Cayley 
graph of an \textbf{amenable}\,\footnote{thus $\GG$ is in $\OOO$} group, then 
every cluster of $\mu_{p}$-a.e.\ subgraph is in $\OOO$}.
\end{quote}

\noindent
The main goal of this paper is to complete these results and prove the following:
\begin{theoremintro}\label{thmintro: harmonic Dirichlet functions
    in uniqueness phase} (Theorem~\ref{th: restriction and clusters in OHD})
 Let $\GG$ be a Cayley graph of a finitely
    generated group.  Consider Bernoul\-li percolation in the uniqueness
    phase.
    Then $\mu_{p}$ a.s.\ the infinite cluster $\omega_{\infty}$ of
    $\omega$ satisfies:\\
$\omega_{\infty}$ has no harmonic Dirichlet functions besides the constants if and
    only if $\GG$ has no harmonic Dirichlet functions besides the constants:
    \textrm{$\mu_{p}$ a.s.\ } $$ \omega_{\infty} \in \OOO
    \Longleftrightarrow \GG\in \OOO$$
\end{theoremintro}
This answers question~1.10 and settles conjecture~1.11 of 
\cite{BLS99}.  Together with their Corollary~4.7 about the nonuniqueness phase, this theorem allows
to complete the picture for Cayley graphs:

\begin{corollaryintro}\label{corintro: uniqueness and nonuniqueness 
phase}
    Let $\GG$ a Cayley graph of a finitely generated group.
   \begin{itemize}
       \item  If $\GG$ is not in $\OOO$ (i.e.\  $\beta_{1}(\Gamma)\not=0$), 
       then a.s.\,\footnote{every infinite cluster for $\mu_{p}$-almost every subgraph $\omega\in 
       \{0,1\}^{\EE}$} the infinite clusters of 
       Bernoulli bond percolation (for both nonuniqueness and 
       uniqueness phases) are not in $\OOO$.
       \item If $\GG$ is in $\OOO$ (i.e.\  $\beta_{1}(\Gamma)=0$), then  
       a.s.\footnotemark[\thefootnote] the infinite clusters of 
       Bernoulli bond  percolation are\\
    -- not in $\OOO$ in the nonuniqueness phase\\
    -- in $\OOO$ in the uniqueness phase.
   \end{itemize}
\end{corollaryintro}
Since percolation at $p_{c}$ belongs to the finite phase for 
nonamenable Cayley 
graphs, it turns out that every $p_{0}>p_{c}$ 
suits in Th.~1.12 of \cite{BLS99} recalled above, while $p_{0}=p_{c}$ 
doesn't.

\medskip
In the course of the proof, a crucial use is made of the notion of
(first) $L^{2}$ Betti numbers for measured equivalence relations,
introduced in \cite{Gab02}. We consider two standard equivalence relations 
with countable classes associated with our situation: the \textbf{full equivalence relation} $\RRfull$
and the \textbf{cluster equivalence relation} $\RRcluster$.
They are defined concretely or also more geometrically (see Section~\ref{subsect: full equiv. 
rel. for Cayley graphs} and \ref{subsect: cluster equiv rel})
in terms
of two \textit{laminated spaces} $\Lfull$ and $\Lcluster$, constructed from
$\Omega\times \GG$ after taking the quotient under the diagonal $\GGamma$-action
and removing certain edges. Their laminated structure comes
from the fact that these spaces are equipped with a measurable partition into
\textit{leaves}, corresponding to the decomposition of $\Omega\times \GG$ into the graphs
$\{\omega\}\times\GG\simeq \GG$.

\medskip
One shows (Section~\ref{sect:fields graphs, HD functions and L2 Betti
numbers}) that the first $L^{2}$ Betti number of such an equivalence
relation, generated by such a
1-dimensional lamination (in fact generated by a \textbf{graphing} in
the sense of \cite{Lev95,Gab00} or Section~\ref{sect: background Measured Equivalence Relations} or
example~\ref{examp: unoriented graphings} below), vanishes if and only
if the leaf of almost every point in the transversal is a graph without harmonic Dirichlet functions, besides the constants:
\begin{theoremintro}\label{thmintro: Betti 1 and graphings} (Corollary~\ref{cor: relate beta-1 and HD functions for relations})
    Let $\RR$ be a measure-preserving equivalence relation on the
    standard Borel probability measure space $(X,\mu)$.  Let $\Phi$ be
    a graphing generating $\RR$.  If the graph $\Phi[x]$ associated
    with $x\in X$ has $\mu$ a.s.\ bounded degree\,\footnote{a bound on the number of neighbors of each vertex}, then
    $\beta_{1}(\RR,\mu)=0$ if and only if $\mu$
    a.s.\  $\HD(\Phi[x])=\Cmath$.
\end{theoremintro}

Now the triviality of the first $L^{2}$ Betti number of an equivalence
relation is invariant when taking a restriction to a Borel subset that
meets almost every equivalence class \cite[Cor.~5.5]{Gab02}.  
Denote by ${\Omegaunique}$ the subset of $\omega$'s with a
\textbf{unique} infinite cluster and such that the base point $\origin$
belongs to that cluster.  In the uniqueness phase, ${\Omegaunique}$
meets almost every $\RRfull$-equivalence class and the restrictions of
$\RRfull$ and $\RRcluster$ to ${\Omegaunique}$ define the same
equivalence relation.
And
Theorem~\ref{thmintro: harmonic Dirichlet functions in uniqueness
phase} follows.

\bigskip
Invariant bond percolation on a locally finite graph $\GG$, for a
group $\HH$ of automorphisms of $\GG$, is also considered  in a more
general setting than just Cayley graphs.  For the invariance property
of the measure to be of any use, the group has to be big enough.  The
standard hypothesis is that $\HH$ is \textbf{transitive} or at least
\textbf{quasi-transitive} (there is only one, resp.\  only
finitely many orbits of vertices).

When closed in all automorphisms of $\GG$, the group $\HH$ is locally
compact and equipped with a unique (up to multiplication by a
constant) left invariant Haar measure.  If that measure is also right
invariant, then $\HH$ is called \textbf{unimodular}. 
A graph with a unimodular quasi-transitive group $\HH$, is called itself
unimodular. The
unimodularity assumption is a quite common hypothesis in
invariant percolation theory, where it is used in order to apply a simple
form of the \textit{mass-transport principle} (see for example
\cite[sect. 3]{BLPS99a} and Section~\ref{sect:Invariance of the Measure and Unimodularity}).

The same unimodularity assumption appears here, for a related reason:
in order to ensure that a certain equivalence relation preserves the
measure (see Section~\ref{sect:Invariance of the Measure and
Unimodularity}).  We obtain the following generalization\,\footnote{Compare with \cite{BLS99} where Theorem 5.7 extends Theorem
1.12 (recalled above) to the more general setting of a unimodular transitive 
graph, and for much more general percolations than Bernoulli percolation.} of
Theorem~\ref{thmintro: harmonic Dirichlet functions in uniqueness
phase}:
\begin{theoremintro}\label{theo intro: when are selectable clusters in OHD} 
    Let $\GG$ be a locally finite graph, $\HH$ a closed
    transitive unimodular group of automorphisms of $\GG$ and $\mathbf{P}$
    any $\HH$-invariant percolation.  On the Borel subset of subgraphs
    with finitely many infinite clusters, $\mathbf{P}$-almost surely the
    infinite clusters belong (resp.\   do not belong) to $\OOO$ if and only
    if $\GG$ belongs (resp.\   doesn't belong) to $\OOO$.
\end{theoremintro}
The clusters of this theorem satisfy the more general property 
(to be introduced in Section~\ref{sect:about selectability}) of being 
\emph{(virtually) selectable} and the proof is given in that
context (Section~\ref{sect: Harmonic
Dirichlet Functions and Clusters for Transitive Graphs},
Theorem~\ref{th: virt select clusters in/out OHD}). 
\bigskip

The most studied invariant percolation, beyond Bernoulli, is probably the \textbf{Random-Cluster Model}. It was introduced by C.~Fortuin
and P.~Kasteleyn \cite{FK72} in relation with Ising and Potts models
 as explained for instance in \cite[Prop.~2.3 and 2.4]{HJL02a}.

It is a (non-independent) percolation process, governed by two 
parameters\,\footnote{The temperature, $T$ in Ising or Potts 
models is linked with the 
parameter $p$ of the Random-Cluster model by $p=1-e^{-\frac{2}{T}}$. 
The parameter $q$, taken to be $q=2$ in the Ising Model, resp.\  $q\in \Nmath$ in the Potts Model, may assume any
value in $[1,\infty)$ for the Random-Cluster model.
For example, the 
(free) Gibbs distribution $\mathtt{FPt}_{\frac{1}{T},q}$ of the Potts model
on $\{0,1\}^{\VV}$ is obtained from $\FRC_{p,q}$ by choosing a subgraph $\omega\in \{0,1\}^{\EE}$
according to $\FRC_{p,q}$ and then choosing a color in 
$\{1,2,\cdots,q\}$ uniformly and independently on the vertices of each cluster.}
$p\in [0,1]$ and $q\in [1,\infty]$. It is defined through a limit procedure by considering an exhaustion $\GG_m$
of $\GG$ by finite subgraphs, and on the set of subgraphs of $\GG_m$, this measure only
 differs from the Bernoulli$(p)$ product measure by the introduction 
of a weight 
($q$ to the power the number of clusters). However, the count of this number 
of clusters
is influenced by the boundary conditions. This leads to 
two particular incarnations of the Random-Cluster model: $\WRC_{p,q}$ and $\FRC_{p,q}$ according to the {\em Wired} (the boundary points are all connected from the exterior)
or {\em Free} (there is no outside connection between the boundary points)
boundary conditions .
These invariant bond percolations both exhibit phase transitions,
for each $q$, similar to that of Bernoulli percolation, leading to critical values $p_c(q)$ and $p_u(q)$ (denoted more precisely by $p_c^{\mathtt{W}}(q)$, $p_u^{\mathtt{W}}(q)$  and $p_c^{\mathtt{F}}(q)$,  $p_u^{\mathtt{F}}(q)$ in case the boundary conditions have to be emphasized). They  ``degenerate'' to Bernoulli percolation when $q=1$.

The reader is invited to consult the papers \cite{HJL02a, HJL02b}
of O.~H\"aggstr\"om, J.~Jonasson and R.~Lyons,
for most of the results relevant for this paper and for details and further references.

\medskip
\noindent
The above Theorem~\ref{theo intro: when are selectable clusters in OHD} 
obviously specializes to:
\begin{corollaryintro}
    Let $\GG$ be a locally finite graph admitting a transitive unimodular
     group of automorphisms.
    Consider the Random-Cluster model $\RC_{p,q}=\WRC_{p,q}$ or $\FRC_{p,q}$
    in the uniqueness phase. 
	Then $\RC_{p,q}$ a.s.\ the infinite cluster admits non- (resp.\  only) constant harmonic 
	Dirichlet functions if and
    only if $\GG$ admits non- (resp.\  only) constant harmonic Dirichlet functions.
\end{corollaryintro}

\subsection{On the nonuniqueness Phase}

One of the most famous conjectures in the subject is probably Conjecture 6
 of Benjamini and Schramm \cite{BS96}:
\begin{quote}
    \emph{The nonuniqueness phase always exists\,\footnote{i.e.\  $p_c(\GG) < p_u(\GG)$ for Bernoulli percolation} for 
Cayley graphs $\GG$ of nonamenable groups.} 
\end{quote}
	and more generally
\begin{quote}
	\emph{If the quasi-transitive graph $\GG$ has positive Cheeger constant, then $p_c(\GG) < p_u(\GG)$.}
\end{quote}

I. Pak and T. Smirnova-Nagnibeda \cite{PSN00} proved that 
each finitely generated nonamenable group admits a Cayley graph for which $p_c<p_u$. On the other hand, the groups with \textit{cost} 
strictly bigger than~1 (see \cite{Gab00} or
Section~\ref{sect: background Measured Equivalence Relations} below, item ``cost'') are the only ones
for which it is known that $p_{c}\not=p_{u}$ for every Cayley graph
(R.~Lyons \cite{Lyo00}). This class of Cayley graphs contains all those outside $\OOO$
(see Th.~\ref{thm: beta-1 and HD of a Cayley graph} and \cite[Cor.~3.23]{Gab02}), but it is unknown whether the reverse inclusion holds.
\medskip

We are able, using our $\ell^2$ methods, to extend Lyons' result to the unimodular setting and to make some progress for Random-Cluster model.
Our treatment doesn't make use of the continuity of the 
probability that $\origin$ belongs to an infinite cluster, 
but only of the expected degree.
We show:
\begin{theoremintro}
{(Cor.~\ref{cor: pc<pu for G not in OHD})}
\label{thintro: pc<pu for G not in OHD} 
     Let $\GG$ be a unimodular transitive locally finite graph.
    If $\GG$ doesn't belong to $\OOO$, then the nonuniqueness phase 
    interval of Bernoulli percolation has non-empty interior:
    $$p_{c}(\GG) < p_{u}(\GG)$$
\end{theoremintro}

\medskip
In fact, to each unimodular transitive locally finite graph $\GG$, 
we associate (see def.~\ref{def: beta1 of the graph}) a \emph {numerical invariant}~$\beta_{1}(\GG)$,
which can be interpreted as the first $\ell^2$ Betti number of any 
closed transitive group of automorphisms of $\GG$. \emph{It vanishes if 
and only if $\GG$ belongs to $\OOO$}. In case $\GG$ is a Cayley graph 
of a group $\Gamma$, then $\beta_1(\GG)=\beta_1(\Gamma)$.
A transitive tree of degree $d$ has  $\beta_1(\GG)=\frac{d}{2}-1$.

\medskip
For Bernoulli percolation, we get the more precise estimate, where $\val(\GG)$ denotes the degree
 of a (any) vertex of $\GG$:
$$0<\beta_{1}(\GG)\leq \frac{1}{2} \val(\GG) (p_{u}(\GG)-p_{c}(\GG)).$$
Observe 
that $\val(\GG)\, p $ is the expected degree 
$\mu_p [\val (\origin)]$ of a base point $\origin$ with respect to the 
Bernoulli measure of parameter $p$.
These inequalities appear as Corollary~\ref{cor: pc<pu for G not in OHD} 
of a quite general result
(Th.~\ref{th: general pc<pu+ couplings}) which applies to more general percolations, like
the free or the wired Random-Cluster model, $\RC_{p,q}=\WRC_{p,q}$ or $\FRC_{p,q}$:
\begin{theoremintro}
{(Cor.~\ref{cor: `pc<pu' for G not in OHD and RC model})}
\label{thintro: beta-1 and pc,pu for Random-Cluster}
	Let $\GG$ be a unimodular transitive locally finite graph, not in $\OOO$.
	Fix the parameter $q\in [1,\infty)$.
	The gap between the left limit (when $p\nearrow p_c(q)$) and the right limit (when $p\searrow 
	p_u(q)$) of the expected degree of a base point $\origin$ with respect to the measure $\RC_{p,q}$ 
	satisfies:
	$$0<\beta_1(\GG)\leq \frac{1}{2} \Big( \RC_{p_u +, q} [\val (\origin)]-\RC_{p_c -, q} [\val 
	(\origin)]\Big).$$
\end{theoremintro}
Indeed the function 
$p\mapsto \FRC_{p,q}[\val (\origin)]:=\int_{\{0,1\}^{\EE}} \val(\origin)(\omega)
d\FRC_{p,q}(\omega)$ is left continuous while $p\mapsto \WRC_{p,q}[\val (\origin)]$ is right continuous in $p$, but the possible remaining discontinuity doesn't allow to conclude in general that $p_c(q)<p_u(q)$.
Of course, $\RC_{p_u +, q} [\val (\origin)]$ denotes $\lim_{p\searrow p_u(q)}\RC_{p,q}[\val(\origin)]$
and 
$\RC_{p_c -, q} [\val 
	(\origin)]:=\lim_{p\nearrow p_c(q)}\RC_{p,q}[\val(\origin)]$.
	
\bigskip
There is no fundamental qualitative difference when shifting 
from transitive to quasi-transitive graphs. The slight modifications
are presented in section~\ref{sect: quasi-transitive graphs}.

\subsection{About Higher Dimensional Invariants and Treeablility}

Higher dimensional $\ell^2$ Betti numbers are also relevant in percolation 
theory.
Y.~Peres and R.~Pemantle \cite{PP00} introduced the 
percolation theoretic notion of countable
\textbf{treeable groups}. They are groups $\Gamma$ for which the space
of trees with vertex set $\Gamma$ admits a $\Gamma$-invariant
probability measure. 
They proved that nonamenable direct products are not treeable.
It is not hard to show that being treeable is equivalent
to being \textit{not anti-treeable} in the sense of \cite[D\'ef.~VI.1]{Gab00} or 
to having \textit{ergodic dimension $1$} in the sense of \cite[D\'ef.~6.4]{Gab02}.

Similarly, R.~Lyons introduced the notion of \textbf{almost treeable groups}:
They are groups for which the space of forests with vertex set $\Gamma$
admits a sequence $(\mathbf{P}_{n})_{n\in \Nmath}$ of $\Gamma$-invariant probability measures with the property that for each pair
of vertices, the probability that they belong to the same connected component
of the forest tends to $1$ as $n$ tends to infinity:
$\forall \gamma_1,\gamma_2\in \Gamma$, $\lim_{n\to \infty} \mathbf{P}_n(\gamma_1\leftrightarrow \gamma_2)=1$.
Clearly treeable implies almost treeable.

It is not hard to show that $\Gamma$ is almost treeable if and only if it
has \textit{approximate ergodic dimension $1$}; where 
the \textbf{approximate ergodic dimension} of $\Gamma$ is the minimum 
of the approximate dimensions of the equivalence relations produced by a free p.m.p.\ action of $\Gamma$ on a standard Borel space (see \cite[D\'ef.~5.15]{Gab02}).

The next theorem follows from 
\cite[Cor.~5.13, Prop.~5.16, Prop.~6.10]{Gab02}
and imposes serious restrictions 
for a group to be treeable or almost treeable. In particular, 
lattices in $\mathrm{SO}(n,1)$ are treeable if and only if 
$n\leq 2$. Also direct products $\Gamma_1\times \Gamma_2$ are not almost treeable as soon as $\Gamma_1$ and $\Gamma_2$ contain a copy of the free group 
$\FF_2$.
This answer questions of R.~Lyons and Y.~Peres (personal communication).
\begin{theoremintro}
If $\Gamma$ is treeable in the sense of \cite{PP00}, then $\beta_1(\Gamma)=0$ if and only if $\Gamma$ is amenable.
If $\Gamma$ is almost treeable, then its higher $\ell^2$ Betti numbers  all vanish:  $\beta_n(\Lambda)=0$ for every $n\geq 2$.
\end{theoremintro}

\bigskip

\noindent{\textbf{Warning for the reader}}\hskip5pt 
Theorem~\ref{thmintro: harmonic Dirichlet functions in uniqueness
phase} is clearly a specialization of Theorem~\ref{theo intro: when
are selectable clusters in OHD}.  However, for the convenience of
the reader mainly interested in Cayley graphs and also to serve as a
warm-up for the more technical general case, we present first a
separate proof of Theorem~\ref{thmintro: harmonic Dirichlet functions
in uniqueness phase} (Section~\ref{sect:Percolation on Cayley
Graphs}, Subsection~\ref{subsubsect: Bernoulli Percolation} and 
Th.~\ref{th: restriction and clusters in OHD}), while Theorem~\ref{theo intro: when
are selectable clusters in OHD} is proved in Section~\ref{selected Clusters and Harmonic Dirichlet Functions}.  A necessary consequence is a certain number of
repetitions. 
Section~\ref{sect: Non uniq. pahse Harm. Dir. Funct} is devoted to the proof of 
Theorems~\ref{thintro: beta-1 and pc,pu for Random-Cluster} and 
\ref{thintro: pc<pu for G not in OHD}.
Some notions related 
to equivalence relations are recalled in Section~\ref{sect: background Measured Equivalence Relations}.
 Theorem~\ref{thmintro: Betti 1 and graphings} is used at several places. It is proved
as Corollary~\ref{cor: relate beta-1 and HD functions for relations}
(see also Remark~\ref{rem: unoriented vs oriented graphings}). But the 
 sections~\ref{sect:Harmonic Dirichlet Functions and L2 Cohomology} and 
\ref{sect:fields graphs, HD functions and L2 Betti numbers} are quite technical, and 
I put it back until the end of the paper. 
It may be a good advice to skip them and 
to keep Theorem~\ref{thm: Betti 1 and graphings} and 
Corollary~\ref{cor: relate beta-1 and HD functions for 
relations} as ``black boxes'' for a first reading.

\bigskip

\noindent{\textbf{Acknowledgment}}\hskip5pt
I wish to thank Thierry Barbot,  \'Etienne Ghys, 
Alekos Kechris, Yuval Peres and Bruno S\'evennec for many valuable discussions,
and especially Vincent Beffara who helped me to understand several references.
I am particularly grateful to Russell Lyons who explained to me many 
entertaining and impressive results in graph percolation theory, as well as some 
connections with my previous work on cost of equivalence relations, and who
suggested that the notion of $\ell^2$ Betti numbers of equivalence 
relations could be relevant for question~1.10 and conjecture~1.11 of \cite{BLS99}. Thanks also for his careful proofreading.

My acknowledgments are due to the anonymous referee who recommended the addition
of a section about the quasi-transitive case (Section~\ref{sect: quasi-transitive graphs} of the current version) and for his accurate and pertinent comments.

\medskip
\noindent To Carole, Valentine, Alice, Lola and C\'elestin.

\tableofcontents

\section{Percolation on Cayley Graphs}\label{sect:Percolation on Cayley Graphs}

Let $\GG=(\VV,\EE)$ be a Cayley graph of a finitely generated group
$\GGamma$. Let $\origin$ be a base vertex, for example the vertex representing the 
identity element of $\GGamma$. The group $\GGamma$ acts freely on the 
set $\EE$ of edges and freely transitively on the set $\VV$ of vertices.

\medskip

The space ${\Omega}:=\{0,1\}^{\EE}$ is the space of \textbf{colorings} 
(assignment of a number) of the edges of $\GG$
with two colors ($0$ and $1$) of $\GG$.
A point $\omega\in \Omega$ is also the characteristic function of a subset of $\EE$. 
When viewed as a subgraph of $\GG$ it is denoted by $\omega(\GG)$. It then 
has the same set of vertices $\VV$
as $\GG$ and for edges the set of \textbf{retained} or \textbf{open} edges:  
those edges $e\in \EE$ with color $\omega(e)=1$.
It has the same base vertex $\origin$ as $\GG$.
The \textbf{cluster} $\omega(v)$ of a 
vertex $v$  is its connected component in $\omega(\GG)$.
The action of $\GGamma$ on $\EE$ induces an action\,\footnote{
For $\gamma\in \GGamma$: $\omega'=\gamma\cdot \omega$ if and only if 
$\omega'(e)=\omega(\gamma^{-1} e)$ for every edge $e\in \EE$.} on the space ${\Omega}$ of colorings.

\medskip
Let $(X,\mu)$ be a standard Borel probability space together with 
\begin{itemize}
    \item a probability measure-preserving (p.m.p.) action of 
    $\GGamma$, which is  \emph{(essentially) free}\,\footnote{the Borel set of points $x\in X$ with non-trivial stabilizer 
    have $\mu$-measure $0$}, and
    \item a $\GGamma$-equivariant Borel map $\pi:X\to \{0,1\}^{\EE}$.
\end{itemize}
The push-forward measure $\pi_{*}\mu$ is a \textbf{$\GGamma$-invariant 
bond percolation} on $\GG$.

\subsection{The Full Equivalence Relation}\label{subsect: full equiv. 
rel. for Cayley graphs}
Consider now the space $X\times \GG$ with the diagonal action of 
$\GGamma$. It is a ``laminated space'', with leaves $\{x\}\times \GG$.

Dividing out by the diagonal action of $\GGamma$, one gets the 
laminated space $\Lfull=\GGamma\backslash (X\times \GG)$: the 
\textbf{full lamination}. 
It is a (huge, highly disconnected) graph with vertex set $\GGamma\backslash 
(X\times \VV)$ and edge set $\GGamma\backslash (X\times \EE)$.
A \textbf{leaf} is a connected component of this graph.

\noindent
-- Because of the freeness of the $\GGamma$-action on $\VV$, the image $X^{\bullet}$
in $\Lfull$ of the space $X\times \{\origin\}$ is an embedding, leading to a 
natural identification of $X$ with $X^{\bullet}$.\\
-- Because of the transitivity on $\VV$ of the $\GGamma$-action, 
$X^{\bullet}$ equals $\GGamma\backslash (X\times \VV)$. 
\begin{displaymath}
\begin{array}{rcl}
 X &\overset{\sim}{\longrightarrow} & X^{\bullet}=\GGamma\backslash (X\times \VV)\\
 x & \longmapsto &
(x,\origin)\sim (\gamma x, \gamma \origin)
\end{array}
\end{displaymath}

Let's denote by 
$\mu^{\bullet}$ the push-forward of the measure $\mu$ to $X^{\bullet}$.

With the (any) choice of $\origin$, the $\GGamma$-set $\VV$ identifies with $\GGamma$ equipped with the action by left multiplication
and 
the left action of $\GGamma$ on itself, by multiplication by the inverse
on the right,
induces on $X^{\bullet}=\GGamma\backslash (X\times \VV)\simeq \GGamma\backslash 
(X\times \GGamma)$ a $\mu^{\bullet}$-preserving $\GGamma$-action, isomorphic to the original one on $X$. 
With $\origin\leftrightarrow 1$ and  $\gamma \origin \leftrightarrow \gamma$ one gets 
$\gamma_{1} (x,1):=(x, 1 \gamma_{1}^{-1})=(x, \gamma_{1}^{-1} 
 1)\sim (\gamma_{1} (x), 1)$.

\noindent
-- Because of the freeness of the $\GGamma$-action on $X$, the leaf of 
$\mu^{\bullet}$-almost every $x^{\bullet} \in  X^{\bullet}$ is isomorphic to 
$\GG$.

\begin{definition}
Define the \textbf{full equivalence relation} $\RRfull$ on $X^{\bullet}$ by 
$x^{\bullet}\RRfull y^{\bullet}$ if and only if $x^{\bullet}$ and 
$y^{\bullet}$ are vertices of the same $\Lfull$-leaf. 
\end{definition}
It is isomorphic
with that given by the $\GGamma$-action on $X$ and preserves the measure.

\medskip
\begin{tabular}{lp{14cm}}
    $\RRfull$: &
 \emph{Two points $x,y$ are $\RRfull$-equivalent if and only if 
 there is $\gamma\in \GGamma$ such that $\gamma x=y$.}
\end{tabular}

\begin{example}\label{examp: Gamma=Z}
Let $\GG=L$ be the standard Cayley graph of $\Zmath$: the simplicial line.
Take $\origin$ as the point $0$ of the line.
The product space $X\times L$ is a kind of cylinder laminated by lines.
A fundamental domain for the diagonal $\Zmath$-action on $X\times L$
is given by $X\times [0,1)$.
Denote by $t$ the automorphism of $X$ given by the generator $1$ of $\Zmath$.
The quotient space $\Zmath\backslash(X\times L)$ identifies with 
the usual \emph{suspension} or \emph{mapping torus} construction
$(\omega,0)\sim(t\omega,1)\backslash (X\times [0,1])$
obtained from $X\times [0,1]$ (laminated by $\{\omega\}\times [0,1]$)
 by gluing together the top and bottom
levels $X\times \{0\}$ and $X\times \{1\}$ after twisting by $t$.
The gluing scar is the transversal~$X^{\bullet}$.
\end{example}

\subsection{The Cluster Equivalence Relation}\label{subsect: cluster equiv rel}

Now, thanks to the map $\pi:X\to \{0,1\}^{\EE}$,  the field of graphs $x\mapsto \{x\}\times \GG$
becomes a $\GGamma$-equivariant field of colored graphs $x\mapsto \pi(x)$. Each leaf of 
$\Lfull$ becomes a colored graph. 

By removing all the $0$-colored edges, one defines a subspace $\Lcluster$ of $\Lfull$: the 
\textbf{cluster lamination}. A leaf of $\Lcluster$ is a connected 
component of $1$-colored (or retained) edges.

\begin{definition}
Define the \textbf{cluster equivalence relation} $\RRcluster$ on $X^{\bullet}$ by 
$x^{\bullet}\RRcluster y^{\bullet}$ if and only if $x^{\bullet}$ and 
$y^{\bullet}$ are vertices of the same $\Lcluster$-leaf. It is a 
subrelation of $\RRfull$.
\end{definition}

For $\mu$-almost every $x\in X$, the leaf of $x^{\bullet}$ is isomorphic to the 
cluster $\GG_{x}:=\pi(x)(\origin)$ of the vertex $\origin$ in the subgraph $\pi(x)$ of $\GG$. 
Thus the $\RRcluster$-class of $x^{\bullet}$ is infinite if and only if 
the corresponding cluster $\pi(x)(\origin)$ is infinite. 
For each $x^{\bullet}\in X^{\bullet}$, the family of $\RRcluster$-classes into which its $\RRfull$-class 
decomposes is in natural bijection with the clusters of $\pi(x)$.
The $\RRfull$-class of $x^{\bullet}$ contains $n$ infinite $\RRcluster$-classes i{f}{f}
$\pi(x)$ has $n$ infinite clusters.

\medskip

Let $e=[\origin,\gamma^{-1} \origin]$ be an edge with end point $\origin$.
Once descended in $\Lfull$, the edge 
$\{x\}\times e=[(x,\origin), \underbrace{(x,\gamma^{-1}
\origin)}_{\sim(\gamma x,\origin)}]$
is retained in $\Lcluster$ i{f}{f} ${\pi(x)(e)=1}$.
In this case, the vertices $\origin,\gamma^{-1}\origin$ are in the same cluster of 
$\pi(x)$, while $x,\gamma x$ are $\RRcluster$-equivalent.
More generally:

\medskip
\begin{tabular}{lp{14cm}}
    $\RRcluster$: &
    \emph{Two points $x,y$ are $\RRcluster$-equivalent if and only if 
    there is $\gamma\in \GGamma$ such that $\gamma x=y$ and the vertices $\origin, \gamma^{-1} \origin$ 
    are in the same cluster of $\pi(x)$.}
\end{tabular}
\medskip

\noindent
It may be relevant to emphasize the role of $\pi$  and
include it in the notation: $\RRcluster_{\pi}$.

\medskip
\noindent
\textbf{Example~\ref{examp: Gamma=Z} (continued)}
\emph{The edges are divided into the positive ones and the negative one according to their position with respect
to $\origin=0$.
The $\RRcluster$-equivalence class of $x\in X$ consists of the iterates $t^{k}(x)$ for 
$k\in\{-i, -i+1,  \cdots, 0, 1, 2, \cdots, j-1, j\}$, where $j$ and $i$ are the number
of negative and positive edges colored $1$ in $\pi(x)$, 
starting from $\origin$:\\
$\pi(x)=
(
\overbrace
{\cdots, 0, \underbrace{1,1,1,\cdots,1}_{j \mathrm{\ times \ }}}
^{\mathrm{\ negative\ edges\ }},
\overbrace{\underbrace{1,1,1,\cdots,1}_{i \mathrm{\ times \ }}, 0, \cdots}^{\mathrm{\ positive\ edges\ }})$.
}

\bigskip
\emph{From now on, and until the end of Section~\ref{sect:Percolation on Cayley Graphs}
we won't distinguish between $(X,\mu)$ and $(X^{\bullet}, \mu^{\bullet})$}.

The \textbf{uniqueness set} $\Omegaunique$ is the Borel subset of points 
$\omega$ of $\Omega$ such that $\omega$ has a \textbf{unique} infinite cluster and such that $\origin$ belongs to it.
The \textbf{uniqueness set $\Xunic$ of $\pi$} is the Borel subset of points of $X$ such that 
$\pi(x)$ belongs to $\Omegaunique$.

\begin{proposition}\label{prop:Rfull-Rcluster coincide on uniqueness set}
    When restricted to the uniqueness set $\Xunic$ the two equivalence relations 
    $\RRfull$ and $\RRcluster$ do coincide.
\end{proposition}
The point here is the \textit{selectability} of the infinite cluster (see Subsection~\ref{sect:about 
selectability} devoted to that subject).

\proof Let $x,y$ be two $\RRfull$-equivalent points of $\Xunic$. 
Since $\pi(x)$ contains only one infinite cluster, the $\RRfull$-class 
of $x$ contains only one infinite $\RRcluster$-class. The
$\RRcluster$-equivalence classes of $x$ and $y$ being both infinite, 
they coincide.
\endofproof

\subsection{Examples}

\subsubsection{Trivial Example}\label{subsubsect: trivial example}
Apply the above construction to the particular constant map
 $x\overset{\pi_{1}}{\mapsto} \GG$, sending every point $x$ to 
the full graph $\GG$ ($\omega\equiv 1$, i.e.\  $\pi(x)(e)=1$ for every $x\in 
X$ and $e\in \EE$, which is fixed by the whole of $\GGamma$).
 In this case, almost every $\GG_{x}$ is 
just $\GG$, the uniqueness set $\Xunic$ equals $X$ and $\RRcluster=\RRfull$.

\subsubsection{Bernoulli Percolation}\label{subsubsect: Bernoulli 
Percolation}
The main example is given by $X=\{0,1\}^{\EE}$ itself, $\pi=id$
and $\mu=\mu_{p}$ the Bernoulli measure with survival parameter $p$,
i.e.\  $\mu_{p}$ is the product of the measures giving weights
$1-p, p$ to $0,1$.
The $\GGamma$-action is essentially free for $p\not= 0$ or $1$.

Stricto sensu, this example is just what is needed for the statement of
Theorem~\ref{thmintro: harmonic Dirichlet functions in uniqueness phase}. 
However, it is useful to introduce more objects in order to 
better distinguish between the various roles played by the space $\Omega$.

The parameter $p$ \textbf{belongs to the uniqueness phase} if and only if 
$\mu_{p}$-almost every graph in $\Omega$ has a unique infinite 
cluster. In this phase, the uniqueness set $\Omegaunique$
has non-zero $\mu_{p}$-measure and is $\mu_{p}$-a.s.\ the 
union of the infinite $\RRcluster$-classes.
The restrictions of $\RRfull$ 
and $\RRcluster$ to $\Omegaunique$ coincide (Proposition~\ref{prop:Rfull-Rcluster coincide on uniqueness set}).
Theorem~\ref{thmintro: harmonic Dirichlet functions in uniqueness phase}
of the introduction will thus be a corollary of Theorem~\ref{th: restriction and clusters in OHD}
for $p<1$ and is trivial for $p=1$.

\subsubsection{Actions Made Free}\label{subsect: actions made free}
If one is considering a percolation for which the $\GGamma$-action is 
not free, one can switch to a free action by 
taking any probability measure-preserving free $\GGamma$-action on a space 
$(Y,\nu)$, and replacing $\Omega$ by its product with $Y$, equipped
with the product measure and the measure-preserving diagonal action of 
$\GGamma$, together with the natural 
$\GGamma$-equivariant projection $\pi : X=\Omega\times Y \to \Omega$.
General (non-free) percolations are thus treated together in the 
same framework.

\subsubsection{Standard Coupling}\label{examp: standard coupling}

The standard coupling is a very useful way to put all the Bernoulli 
measures $\mu_{p}$ together and 
to vary the map $\pi$ instead of changing the measure on $\Omega$ 
with the parameter $p$.

Let $X=[0,1]^{\EE}$ be the product space with the 
product measure $\mu$ of Lebesgue measures on the intervals $[0,1]$. 
An element of $X$ gives a \textbf{colored graph}: a coloring of the graph $\GG$, 
with  $[0,1]$ as set of colors.
For each $p\in [0,1]$, let $\pi_{p}$ be the $\GGamma$-equivariant map
sending $[0,1]$-colored graphs to $\{0,1\}$-colored ones by
retaining only the edges colored in $[0,p]$:
$$
\pi_{p}:\left( 
\begin{array}{rcl}
    [0,1]^{\EE} & \longrightarrow & \{0,1\}^{\EE} \\
    x   & \longmapsto & \pi_{p}(x) : 
\begin{cases}
\pi_{p}(x)(e)=1 & \text{if $x(e)\in [0,p]$}\\
\pi_{p}(x)(e)=0 & \text{if $x(e)\in (p,1]$}
\end{cases}
\end{array}\right)
$$
Clearly, $\pi_{p}$ pushes the measure $\mu$ to $\mu_{p}$. For each 
value of $p$, one gets the cluster equivalence relation 
$\RRcluster_{p}$, also defined as follows:

\medskip
\begin{tabular}{lp{14cm}}
    $\RRcluster_{p}$: &
    \emph{Two $[0,1]$-colored graphs $x,y$ are $\RRcluster_{p}$-equivalent if and only if 
    there is $\gamma\in \GGamma$ such that $\gamma x=y$ and the vertices $\origin, \gamma^{-1} \origin$ 
    are connected in the colored graph $x$ by a path of edges with colors $\leq p$.}
\end{tabular}
\medskip

    This gives an intuitive picture of the clusters evolution 
as $p$ varies:
The family $(\RRcluster_{p})_{p\in [0,1]}$ is 
strictly increasing. Moreover for every $p$, $\RRcluster_{p}=\cup_{t<p}\RRcluster_{t}$ and  $\RRcluster_{1}=\RRfull$.
The critical value $p_{c}$ is characterized as the supremum of those 
$p$ for which the $\RRcluster_{p}$-classes are finite ($\mu$-a.s.), as 
well as the infimum of $p$ such that $\RRcluster_{p}$ admits a 
$\mu$-non-null set of points with infinite classes.
Much less obvious is the similar characterization of $p_{u}$,
obtained by O.~H\"aggstr\"om and Y.~Peres, who showed that after 
$p_{c}$, \textit{there 
is no spontaneous generation of infinite clusters; all infinite 
clusters are born simultaneously}: 
If $p_c< p \leq q$, then $\mu$-a.s.\  every 
infinite $\RRcluster_{q}$-class contains an infinite $\RRcluster_{p}$-class
\cite{HP99}.
This explains that the uniqueness phase is an interval.
$$\begin{array}{rrl}
    p_{c}&:=&\inf\{p: \hbox{there is a unique infinite cluster for } 
\mu_{p}\}\\
  &= &\sup \{p: \hbox{there is not a unique infinite cluster for } 
\mu_{p}\}
\end{array}$$

\subsubsection{Site Percolation}
An invariant site percolation on $\GG$ is a probability measure 
$\mathbf{P}$ on the space $\{0,1\}^{\VV}$ that is invariant under a 
certain group of automorphisms of $\GG$. To a site percolation 
corresponds a bond percolation by the equivariant map
$\pi:\{0,1\}^{\VV}\to \{0,1\}^{\EE}$ sending a coloring of the 
vertices $\VV$ to the coloring of the edges $\EE$ where an edge 
gets color $1$ if and only if both its endpoints are colored $1$.

\subsubsection{Graphings}\label{subsubsect:examples-graphings}

Let $(\gamma_{1},\gamma_{2}, \cdots, \gamma_{n})$ be the generating system 
defining the Cayley graph $\GG$ and $e_{i}$ be the edge~$[\origin,\gamma_{i}\origin]$.

\medskip
If ${\pi(x)(e_{i})=1}$, the vertices $\origin,\gamma_{i}\origin$ are in the same
cluster of $\pi(x)$, and $x,\gamma_{i}^{-1} x$ are
$\RRcluster$-equivalent.  Define the Borel set $A_{i}:=\{x\in X :
\pi(x)(e_{i})=1 \}$ and the partial Borel isomorphism
$\smash{\varphi_{i}= \gamma_{i\vert A_{i}}^{-1}}$, the restriction of
$\gamma_{i}^{-1}$ to $A_{i}$.  The family ${\Phi=(\varphi_{1},
\varphi_{2}, \cdots, \varphi_{n})}$ is a \emph{graphing} (in the sense
of \cite{Lev95,Gab00} -- see also Section~\ref{sect: background Measured Equivalence 
Relations}) that generates $\RRcluster$: the latter is the smallest
equivalence relation such that $x\sim \varphi_{i}(x)$, for every $x\in
A_{i}$.  For instance, in the above standard coupling (ex.~\ref{examp:
standard coupling}), the cluster equivalence relation $\RRcluster_{p}$
is generated by the graphing $\Phi_{p}=(\varphi_{1}^{p},
\varphi_{2}^{p}, \cdots, \varphi_{n}^{p})$, where $\varphi_{i}^{p}$ is
the restriction of $\gamma_{i}^{-1}$ to $A_{i}^{p}:=\{x\in
[0,1]^{\EE}: x(e_{i})\leq p\}$.

\medskip 
Conversely, given a free p.m.p.\ $\GGamma$-action on
$(X,\mu)$, consider $n$ Borel subsets $A_{i}$, partial isomorphisms
$\smash{\varphi_{i}= \gamma_{i\vert A_{i}}^{-1}}$, the graphing
${\Phi=(\varphi_{1}, \varphi_{2}, \cdots, \varphi_{n})}$ and the
generated equivalence relation $\RR_{\Phi}$.  The coloring
${\pi(x)(e_{i})=1}$ i{f}{f} $x$ belongs to $A_{i}$ extends by
$\GGamma$-equivariance to a map $\pi:X\to \{0,1\}^{\EE}$ whose cluster
equivalence relation $\RRcluster$ coincides with $\RR_{\Phi}$.

\subsection{Harmonic Dirichlet Functions and Clusters for Cayley 
Graphs}\label{subsect: Harmonic Dirichlet Functions and Clusters for Cayley 
Graphs}

We are now able to state the main result of this section.
Recall that we have at hand: (1) A locally finite graph $\GG=(\VV,\EE)$ with a free action of a countable group $\GGamma$, transitive on $\VV$;
(2) A standard probability measure space $(X,\mu)$  with a free 
measure-preserving $\GGamma$-action, and with a 
$\GGamma$-equivariant map $\pi:X\to \{0,1\}^{\EE}$;
(3) The two associated cluster ($\RRcluster$) and full ($\RRfull$) equivalence relations.

\begin{theorem}\label{th: restriction and clusters in OHD}
    For $\mu$-a.e.\ $x$ in the uniqueness set $\Xunic$,
the cluster $\GG_{x}$ of the vertex $\origin$ 
in the colored graph $\pi(x)$ belongs to $\OOO$ if and only if $\GG$ belongs to $\OOO$.
\end{theorem}

For the purpose of proving this result, very little has to be known about the
$L^2$ Betti numbers of equivalence relations. Just assume the 
following ``black box'', which will be further developed in 
Section~\ref{sect:fields graphs, HD functions and L2 Betti numbers}.
The reader feeling more comfortable with the notion of \emph{cost} may 
think at first glance that $\beta_{1}(\RR)=\mathrm{cost}(\RR)-1$
(see Section~\ref{sect: background Measured Equivalence Relations}, item ``cost'').

\begin{enumerate}
    \item[Fact 1.]
    For each measurably defined subrelation $\RR$ of $\RRfull$ on a 
    non-null Borel subset $Y$ of $X$, there is a well-defined 
    notion of first $L^2$ Betti number $\beta_{1}(\RR, \mu_{Y})$, where 
    $\mu_{Y}$ denotes the normalized restricted measure $\frac{\mu_{\vert 
    Y}}{\mu(Y)}$ (\cite{Gab02}).
    In particular, $\beta_{1}(\RRfull, \mu)$, 
    $\beta_{1}(\RRfull_{\vert Y}, 
    \mu_{Y})$ and $\beta_{1}(\RRcluster_{\vert Y}, \mu_{Y})$ 
    are well defined.
    \item [Fact 2.]
    If $Y$ meets almost all $\RRfull$-classes, then  $\beta_{1}(\RRfull, 
    \mu)=\mu(Y) \beta_{1}(\RRfull_{\vert Y}, \mu_{Y})$ \cite[Cor.~5.5]{Gab02}.
    \item[Fact 3.]
    The first $L^2$ Betti number $\beta_{1}(\RRcluster_{\vert 
    Y},\mu_{Y})$ of $\RRcluster_{\vert Y}$ vanishes if and only if for 
    $\mu$-almost every $y\in Y$, the graph $\GG_{y}$ belongs to 
    $\OOO$ (Theorem~\ref{thmintro: Betti 1 and graphings}).
\end{enumerate}

\medskip
\noindent{\textbf{Remark.}}
However, when $Y$ is \emph{$\RRcluster$-saturated} (the $\RRcluster$-class of every 
$y\in Y$ is entirely contained in $Y$), these numbers are ``easily explicitly defined'':
consider the space 
$\HD(\GG_{x})$ of harmonic Dirichlet functions on $\GG_{x}:=\pi(x)(\origin)$.
Its image $\dzero \HD(\GG_{x})$ by the coboundary operator 
$\dzero$ in the $\ell^{2}$ cochains $C_{(2)}^{1}(\GG_{x})$
is a closed subspace of $C_{(2)}^{1}(\GG)$ (every edge outside 
$\GG_{x}$ is orthogonal to it), isomorphic to $\HD(\GG_{x})/\Cmath$.
Denote by $p_{x}:C_{(2)}^{1}(\GG)\to\dzero \HD(\GG_{x})$ 
the orthogonal projection and, for each edge $e\in \EE$, denote by $\mathbf{1}_{e}$ the characteristic function of the 
edge $e$. Let $e_{1}, e_{2}, \cdots, e_{n}$ be a set of 
orbit representatives for the $\GGamma$-action on $\EE$.
\begin{proposition} Let $Y$ be a non-null \emph{$\RRcluster$-saturated} Borel subset of $X$.
    The first $L^{2}$ Betti number of the restricted equivalence relation
    $\RRcluster_{\vert Y}$ on $(Y,\mu_{Y})$ equals
    $$\beta_{1}(\RRcluster_{\vert Y},\mu_{Y})=\frac{1}{\mu(Y)}\sum_{i=1}^{n} 
    \int_{Y}\scal{p_{y}(\mathbf{1}_{e_{i}})}{\mathbf{1}_{e_{i}}} d\mu(y).$$
\end{proposition}
To prove this, we essentially use Theorem~\ref{thm: Betti 1 and graphings} stating 
that $\beta_{1}(\RRcluster_{\vert Y},\mu_{Y})=
\dim_{\RRcluster_{\vert Y}}\int_{Y} \dzero(\HD(\GG_{y})) \ d\mu_{Y}(y)$
and then the definition of the dimension (see \cite[Prop.~3.2 (2)]{Gab02}).
See also Proposition~\ref{prop: computation of L2 Betti numbers}.

\bigskip

\proof (of Theorem~\ref{th: restriction and clusters in OHD})
If $\Xunic$ is a $\mu$-null set, the theorem is empty.
Up to replacing $X$ by the union of the $\GGamma$-orbits meeting 
$\Xunic$, one may assume that $\Xunic$ meets every 
$\RRfull$-class of $X$. The following are then equivalent:
\begin{enumerate}
    \item  $\GG$ is in $\OOO$ \label{enum: G dans OHD}

    \item  $\beta_{1}(\RRfull,\mu)=0$ \label{enum: beta1 de Sfull 
    vanishes}
    
    \item  $\beta_{1}(\RRfull_{\vert \Xunic},\mu_{\Xunic})=0$ \label{enum: beta1 de Sfull IY vanishes}
    
    \item  $\beta_{1}(\RRcluster_{\vert \Xunic},\mu_{\Xunic})=0$ \label{enum: beta1 
    de Scluster vanishes}

    \item  for $\mu$-almost every $x\in \Xunic$, the graph $\GG_{x}$ is in $\OOO$
    \label{enum: a.s. Sigmacluster in OHD}
\end{enumerate}

The  equivalence \ref{enum: G dans OHD}$\iff$\ref{enum: beta1 de Sfull 
vanishes} follows from fact~3 (i.e.\  Theorem~\ref{thmintro: Betti 1 and graphings}) applied to 
the map $\pi_{1}$ of the example~\ref{subsubsect: trivial example}, 
since in this case $X=Y$, $\RRcluster=\RRfull$ and almost every $\GG_{x}$ equals $\GG$.

The equivalence \ref{enum: beta1 de Sfull 
vanishes}$\iff$\ref{enum: beta1 de Sfull IY vanishes} follows from
fact 2.

By Proposition~\ref{prop:Rfull-Rcluster coincide on uniqueness set}, 
$\RRcluster_{\vert \Xunic}=\RRfull_{\vert \Xunic}$.
The key point of the proof is that from \cite{Gab02} these numbers depend only 
on the equivalence relation: one gets  \ref{enum: beta1 de Sfull IY vanishes}$\iff$\ref{enum: beta1 
de Scluster vanishes}.
Again, fact~3 shows the equivalence \ref{enum: beta1 
de Scluster vanishes}$\iff$\ref{enum: a.s. Sigmacluster in OHD},
after noticing that $\mu$ and the normalized measure $\mu_{\Xunic}$ are 
equivalent on $\Xunic$.

It remains to move the quantifier ($\mu$-almost every $x\in \Xunic$)
outside the equivalence \ref{enum: G dans OHD}$\iff$\ref{enum: a.s. Sigmacluster in OHD}.
Let $Y\subset \Xunic$ be the Borel subset of points such that 
the graph $\GG_{y}$ is in $\OOO$. If $Y$ is non-null,
then the argument applied to $Y$ shows that $\GG$ belongs to $\OOO$ and thus 
$Y=\Xunic$ a.s.\ This implies that in case $\GG$ does not belong to 
$\OOO$, then for $\mu$-almost every $x\in \Xunic$, the graph $\GG_{x}$ is 
not in $\OOO$. \endofproof

\begin{remark}
    Observe that the freeness of the $\GGamma$-action on $X$ is a 
    hypothesis made to simplify some arguments ($\gamma$ is the unique element 
    of the group sending $x$ to $\gamma x$) and to apply more directly 
    results from \cite{Gab02}. However, thanks to the example~\ref{subsect: actions made free}, 
    the above Theorem~\ref{th: restriction and clusters in OHD} 
    admits a natural generalization without it.
\end{remark}

\section{Percolation on Transitive Graphs}\label{sect:Percolation on transitive graphs}

Let $\GG=(\VV,\EE)$ be a locally finite transitive\,\footnote{The 
quasi-transitive case is very similar and 
we restrict our attention to the transitive one only to avoid an excess of technicality. The modifications for quasi-transitivity are presented in section~\ref{sect: quasi-transitive graphs}} graph. Let $\mathrm{Aut}(\GG)$ be the automorphism group of 
$\GG$ with the topology of pointwise convergence. Let $\HH$ be a closed
subgroup of $\mathrm{Aut}(\GG)$. We assume that $\HH$ acts 
transitively on the set $\VV$ of vertices. It is locally compact and the stabilizer of each vertex is compact. Let $\origin$ be a base vertex
and denote by $\KK$ its stabilizer.

The action of $\HH$ on $\EE$ induces an action on the space 
${\Omega}=\{0,1\}^{\EE}$ of colorings: for each $\hh\in \HH$, 
$\omega'=\hh\cdot \omega$ if and only if 
$\omega'(e)=\omega(\hh^{-1} e)$ for every edge $e\in \EE$.

\medskip
Let $(X,\mu)$ be a standard Borel probability space together with 
\begin{itemize}
    \item a probability measure-preserving (p.m.p.) action of 
    $\HH$, which is  \emph{essentially free}\,\footnote{the Borel set of points $x\in X$ with non-trivial stabilizer 
    have $\mu$-measure $0$}, and
    \item an $\HH$-equivariant Borel map $\pi:X\to \{0,1\}^{\EE}$.
\end{itemize}

The push-forward measure $\pi_{*}\mu$ is an $\HH$-invariant 
bond percolation on $\GG$.

\subsection{The Full Equivalence Relation}\label{subsect: full equi 
rel for transitive graph}
Consider now the space $X\times \GG$ with the diagonal action of 
$\HH$. It is an $\HH$-equivariant field of graphs 
above $X$, all isomorphic to $\GG$: $x\mapsto \{x\}\times \GG$. It is 
also a ``laminated space'', with leaves $\{x\}\times \GG$.

Dividing out by the diagonal action of $\HH$, one gets the 
laminated space $\Lfull=\HH\backslash (X\times \GG)$: the 
\textbf{full lamination}. 
It is a (huge, highly disconnected) graph with vertex set $\HH\backslash 
(X\times \VV)$ and edge set $\HH\backslash (X\times \EE)$.
A \textbf{leaf} is a connected component of this graph.

Denote by  $X^{\bullet}$ the image
in $\Lfull$ of the space $X\times \{\origin\}$. Because of the transitivity on $\VV$ of the $\HH$-action, $X\times 
\{\origin\}$ meets every $\HH$-orbit of $X\times \VV$, so that
$X^{\bullet}$ equals $\HH\backslash (X\times \VV)$:
\begin{eqnarray}\label{eq: map from (X times *) to X bullet}
\left( \begin{array}{rcl}
 X\simeq X\times \{\origin\} &{\longrightarrow} &X^{\bullet}=\HH\backslash (X\times \VV)\\
 x & \longmapsto &
(x,\origin)\sim (\hh x, \hh \origin)
\end{array} \right) 
\end{eqnarray}
In particular, two points of $X\times\{\origin\}$ happen to be 
identified in $X^{\bullet}$, {i.e.\ } $(x,\origin)\sim (\hh x,\origin)$, if and 
only if $\hh$ belongs to $\KK$. Thanks to the compactness of the 
stabilizer $\KK$ of $\origin$, the space $X^{\bullet}$ gets naturally the 
structure of a standard Borel space (see Proposition~\ref{prop:quotient 
is standard}, Section~\ref{sect:Invariance of the Measure and Unimodularity}).
$$\left( 
\begin{array}{rcl}
   X^{\bullet}=\HH\backslash (X\times \VV) & \simeq &\KK\backslash X\\
(hx,h\gamma \origin)&\mapsto &\KK \gamma^{-1} x
\end{array}
\right) 
$$
Denote by  $\mu^{\bullet}$ the push-forward of the measure $\mu$ to $X^{\bullet}$.
Because of the freeness of the $\HH$-action on $X$, the leaf of 
$\mu^{\bullet}$-almost every $x^{\bullet} \in  X^{\bullet}$ is isomorphic to 
$\GG$.

\medskip
Define the \textbf{full equivalence relation} $\RRfull$ on $X^{\bullet}$ by 
$x^{\bullet}\RRfull y^{\bullet}$ if and only if $x^{\bullet}$ and 
$y^{\bullet}$ are vertices of the same $\Lfull$-leaf.

\medskip
\begin{tabular}{lp{14cm}}
    $\RRfull$: &
    \emph{Two points $x^{\bullet},y^{\bullet}$ are $\RRfull$-equivalent if 
    and only if they admit two representatives in $X\times \VV$ with 
    the same first coordinate, i{f}{f}
    they admit two representatives $x,y$ in $X$ for which there exists 
    $\hh\in \HH$ such that $\hh x=y$, i{f}{f} 
    any of their representatives are in the same $\HH$-orbit.}
\end{tabular}
\medskip

It inherits naturally an \textit{unoriented graphing}
and a \textit{smooth field of graphs}
(see  Section~\ref{sect: background Measured Equivalence 
    Relations}, and examples~\ref{examp: unoriented graphings}, 
\ref{examp: fields of graphs for transitive} of 
Section~\ref{sect:fields graphs, HD functions and L2 Betti numbers})
from the edge set $\HH\backslash (X\times \EE)$, where the graph 
associated with each point admits an 
isomorphism with $\GG$, canonical up to ``rotation
around $\origin$'', i.e.\  up to the action of the stabilizer $\KK$ of $\origin$.

\begin{theorem}\label{th: Rful preserves the measure iff unim.}
    The equivalence relation $\RRfull$ preserves the measure 
    $\mu^{\bullet}$ if and only if the group $\HH$ is unimodular.
\end{theorem}
This result is just an application of 
Theorem~\ref{th: measure-preserved iff unimodular} below. It
sheds another light on the \emph{unimodularity assumption} and on the \emph{Mass Transport Principle} (see the proof of Theorem~\ref{th: measure-preserved iff unimodular}).

\begin{example}
    The simplest example of the graph $\GG$ made of a single infinite line
    is quite eloquent, with $\HH=\mathrm{Aut}(\GG)\simeq\Zmath/2 \Zmath
    \ltimes \Zmath$.  The compact subgroup $\KK=\Zmath/2 \Zmath$ is finite
    and $X^{\bullet}=(\Zmath/2 \Zmath)\backslash X$. If the 
    $\HH$-action on $X$ is ergodic, then the lamination
    $\Lfull$ is not orientable, so that the associated unoriented graphing
    cannot be made (measurably) oriented.
\end{example}

\subsection{The Cluster Equivalence Relation}
\label{subsect:trans graph+ cluster equiv. rel}

Now, thanks to the map $\pi:X\to \{0,1\}^{\EE}$,  the field of graphs $x\mapsto \{x\}\times \GG$
becomes an $\HH$-equivariant field of colored graphs $x\mapsto \pi(x)$ 
so that each leaf of $\Lfull$ becomes a colored graph. 

By removing all the $0$-colored edges, one defines a subspace $\Lcluster$ of $\Lfull$: the 
\textbf{cluster lamination}. A leaf of $\Lcluster$ is a connected 
component of $1$-colored (or retained) edges.

\medskip
Define the \textbf{cluster equivalence relation} $\RRcluster$ on $X^{\bullet}$ by 
$x^{\bullet}\RRcluster y^{\bullet}$ if and only if $x^{\bullet}$ and 
$y^{\bullet}$ are vertices of the same $\Lcluster$-leaf. It is a 
subrelation of $\RRfull$.

The leaf of $\mu^{\bullet}$-almost every $x^{\bullet}$ is a graph 
$\Lfull_{x^{\bullet}}$ which admits an 
isomorphism with the cluster $\GG_{x}:=\pi(x)(\origin)$ of the vertex $\origin$ in the subgraph $\pi(x)$ of $\GG$
for any representative $x$ of $x^{\bullet}$ (just observe that since 
$\kk \origin=\origin$, for any $\kk\in \KK$, the clusters of $\origin$ for $x$
on the one hand and for $\kk x$ on the other hand are isomorphic: $\pi(\kk x)(\origin)=\kk \pi(x)(\origin)$).
Thus the $\RRcluster$-class of $x^{\bullet}$ is infinite if and only if 
the corresponding clusters $\pi(x)(\origin)$ are infinite. 
For each $x^{\bullet}\in X^{\bullet}$, the family of $\RRcluster$-classes into which its $\RRfull$-class 
decomposes is in bijection with the clusters of $\pi(x)$.
The $\RRfull$-class of $x^{\bullet}$ contains $n$ infinite $\RRcluster$-classes i{f}{f}
$\pi(x)$ has $n$ infinite clusters.

\medskip
\begin{tabular}{lp{14cm}}
    $\RRcluster$: &
    \emph{Two points $x^{\bullet},y^{\bullet}$ are $\RRcluster$-equivalent if and only if 
    they admit two representatives in $X\times \VV$ with 
    the same first coordinate $x$ and second coordinates in the same 
    connected component of $\pi(x)$, i{f}{f}
    they admit two representatives $x,y$ in $X$ for which there exists 
    $\hh\in \HH$ such that $\hh x=y$  and the vertices $\origin, \hh^{-1} \origin$ 
    are in the same cluster of $\pi(x)$.}
\end{tabular}
\medskip

Let's check by hand that the above 
characterization doesn't depend on the choice of representatives.
This $h$ defines an isomorphism between the cluster 
$\pi(x)(\origin)=\pi(x)(\hh^{-1}\origin)$ and $\pi(\hh x)(\hh \hh^{-1}\origin)=\pi(y)(\origin)$.
If $\kk_{1} x$ and $\kk_{2} y$ are two other representatives,
$\kk_{1}, \kk_{2}\in \KK$, then $\kk_{2}\hh \kk_{1}^{-1} (\kk_{1} 
x)=\kk_{2} y$ (i.e.\  $h$ has to be replaced by $\kk_{2}\hh \kk_{1}^{-1}$), 
then $\kk_{1} \origin=\origin$ and $\kk_{1} \hh^{-1} \origin= \kk_{1} 
\hh^{-1} \kk_{2}^{-1} \origin= (\kk_{2} \hh \kk_{1}^{-1})^{-1} \origin$ are in the same 
cluster of $\pi (\kk_{1}(x))$. 

\bigskip
The equivalence relation $\RRcluster$ inherits naturally an \textit{unoriented graphing}
and a \textit{smooth field of graphs}
(see  Section~\ref{sect: background Measured Equivalence 
    Relations}, and examples~\ref{examp: unoriented graphings}, 
\ref{examp: fields of graphs for transitive} of 
Section~\ref{sect:fields graphs, HD functions and L2 Betti numbers}).

\begin{remark}
Let $Y^{\bullet}\subset X^{\bullet}$ be the union of the infinite $\RRcluster$-classes.
Assume the action of $\HH$ on $X$ is ergodic.
Then the invariant percolation $\pi_{*}\mu$ has indistinguishable infinite clusters
in the sense of \cite[Sect. 3]{LS99} if and only if the restriction
$\RRcluster_{\vert Y}$ is ergodic.
\end{remark}

\subsection{Measure Invariance, Unimodularity and the Mass-Transport Principle}
\label{sect:Invariance of the Measure and Unimodularity}

Recall that a locally compact second countable group $\GGG$ admits a 
left-invariant Radon measure, its \textit{Haar measure} $\Haar$, unique up to a
multiplicative constant.  Pushed forward by
right-multiplication, the measure is again left-invariant,
and thus proportional to 
$\Haar$. One gets a homomorphism $\mathrm{mod}:\GGG\to \Rmath^{*}_{+}$, 
the \textbf{modular map}, which encodes the defect for $\Haar$ to 
be also right-invariant.
In case the modular map is trivial ($\mathrm{mod}(\GGG)=\{1\}$),
i.e.\  $\Haar$ is also right-invariant, then the group $\GGG$ is called 
\textbf{unimodular}.

\bigskip

Let $(X,\mu)$ be a standard Borel space with a probability measure 
and an essentially free measure-preserving action of a locally compact second 
countable group $\GGG$. Let $K$ be a compact open subgroup of $\GGG$.
Restricted to $K$, the modular function is trivial.

\begin{proposition}\label{prop:quotient is standard}
    The space $\Xbar=K\backslash X$ is a standard 
    Borel space. The quotient map $X\to K\backslash X$ admits a Borel section.
\end{proposition}
\proof The following argument has been explained to me by A.~Kechris.
Pushing forward the normalized  Haar measure $m$ on $K$ by the Borel map (for 
any fixed $x$) $K\to X,\ k\mapsto kx$ defines a measure $m_{x}$ and thus a Borel
map from $X$ to the standard Borel space of probability measures on $X$
(\cite[Th.~17.25]{Kec95}).
But the right invariance of $m$ on $K$  shows\,\footnote{$m_{hx}(A)=m(\{k:khx\in 
A\})=m(\{k'h^{-1}: k'x\in A\})=m(\{k: kx\in A\})=m_{x}(A)$} that 
$x$ and $y$ are in the same $K$-orbit i{f}{f} $m_{x}=m_{y}$.
The $K$-action is then smooth. It follows from \cite[Ex. 18.20]{Kec95}
that the action has a Borel selector.
\endofproof

\bigskip

Let $\RRbar$ be the \textbf{reduced} equivalence relation defined on $\Xbar$ by
$\overline{x}\RRbar \overline{y}$ if{f} $\overline{x}$ and 
$\overline{y}$ admit $\GGG$-equivalent preimages.
Let's denote by $\mubar$ the push-forward probability measure on 
$\Xbar=K\backslash X$, or by $\mubar_{K}$ if one wants to emphasize the 
choice of $K$.  Section~\ref{sect: background Measured Equivalence Relations} 
recalls the terminology for the next Theorem.

\begin{theorem}
	\label{th: measure-preserved iff unimodular}
	The equivalence relation $\RRbar$ on $\Xbar$ is standard countable. 
	It preserves the measure $\mubar_{K}$ if and only if $\GGG$ is unimodular.
\end{theorem}
\proof
It is obviously a Borel subset of $\Xbar\times \Xbar$.
The countability of the classes comes from that of the set 
$K\backslash G$. 

\bigskip
The statement about unimodularity is quite natural 
once one realizes that the decomposition of $\mu$ relatively to 
$\mubar$ makes use of
the right invariant Haar measure on $\GGG$. However, we will follow 
 elementary but enlightening facts leading by two ways to the result. 

Recall (see Section~\ref{sect: background Measured Equivalence Relations})
that $\RRbar$ preserves $\mubar_{K}$ if and only if
the measures $\nu_{1}$ and $\nu_{2}$ on the set 
    $\RRbar\subset X\times X$ coincide, 
    defined with respect to the projections on the first (resp. \ second) coordinate 
    $pr_{1}$ (resp.\  $pr_{2}$) by 
    $\nu_{1}(C)=\int_{X} \#(C\cap pr_{1}^{-1}(x)) d\mubar_{K}(x)$
    and 
    $\nu_{2}(C)=\int_{X} \#(C\cap pr_{2}^{-1}(y)) d\mubar_{K}(y)$.

\bigskip
\textit{Proof by hand}: If $K$, $K'$
are compact open subgroups of $\GGG$, then $K$ is made of
unimodular elements of $\GGG$.  The intersection $K\cap K'$ is a
compact open subgroup of $\GGG$.  By a covering argument, 
its index in $K$ is finite:
$[K:K\cap K']:=\# K/(K\cap K') = \frac{\Haar (K)}{\Haar (K\cap K')}$.
For $\gamma\in \GGG$, $\Haar (\gamma^{-1} K
\gamma)=\mathrm{mod}(\gamma) \Haar (K)$.  In particular, $G$ is
unimodular i{f}{f} all the conjugates of $K$ have the same Haar
measure. Observe that $[K:\gamma K \gamma^{-1}\cap
K]=\mathrm{mod} (\gamma) [K:K\cap \gamma^{-1} K \gamma]$ and that
$K$ and $K'$ as well as their Haar measures are commensurable, so that
the modular function on $\GGG$ is rational.  The reduction map
$(K\cap K')\backslash X\to K\backslash X$ is a.s.\ $[K: K\cap
K']$-to-one, and yields a disintegration of the push-forward measure
$\mubar_{K\cap K'}$ with respect to $\mubar_{K}$, with normalized
counting measure in the fibers.\label{enum: desintegration with
finite fibers}

\noindent
For $\gamma\in \GGG$, consider the graph $C_{\gamma}:=\{(x,\gamma x):x\in X\}$
and its image  in $\RRbar$:
$$\Cbar_{\gamma}:=\{(\overline{x},\overline{\gamma x}):x\in X\}
=\{([K.x],[K.\gamma x]):x\in X\}.$$
For every $k\in K$, $\Cbar_{\gamma}=\Cbar_{k \gamma}$.
Two points $(x_{1},y_{1})$ and $(x_{2},y_{2})$
of the ``curve'' $C_{\gamma}$ define the same point in 
$\Cbar_{\gamma}$
\emph{if{f}} there exist $k,k'\in K$ such that 
$x_{1}=k x_{2}$ and $\gamma x_{1}=y_{1}= k' y_{2}=k' \gamma x_{2}$,
\emph{if{f}}  (since by freeness of the action $k=\gamma^{-1} k' \gamma$) there exists $k\in K\cap (\gamma^{-1} K \gamma)$ such that 
$x_{1}=k x_{2}$, $\gamma x_{1}=y_{1}$ and $\gamma x_{2}=y_{2}$
\emph{if{f}} there exists $k'\in (\gamma K \gamma^{-1}) \cap K$ such that 
$x_{1}=\gamma^{-1} y_{1}$, $x_{2}=\gamma^{-1} y_{2}$ and $y_{1}=k' 
y_{2}$. Thus, the preimage in $C_{\gamma}$ (not in $X\times X$ !) of a point 
in $\Cbar_{\gamma}$ 
is of the form (for certain $x$ and $y$): {$\{(kx,\gamma k x): k\in (K\cap 
\gamma^{-1} K \gamma)\}=\{\gamma^{-1} k' y, k' y): k'\in (\gamma 
K \gamma^{-1}\cap K)\}$}, so that
\begin{center}
$\Cbar_{\gamma}\simeq (K\cap \gamma^{-1} K \gamma)\backslash X$ 
\hfil or \hfil 
$\Cbar_{\gamma}\simeq (\gamma K \gamma^{-1}\cap K)\backslash X$,
\end{center}
according to whether $\Cbar_{\gamma}$ is parameterized by its first 
coordinate or its second coordinate. 
We have thus proved the following:
\begin{proposition}
The projection of $\Cbar_{\gamma}\subset \RRbar$ to the first 
coordinate is $[K:K\cap \gamma^{-1} K \gamma]$-to-one.
To the second coordinate, it is $[K:\gamma K \gamma^{-1}\cap K]$-to-one.
\end{proposition}
It follows that the measures $\nu_{1}$ and $\nu_{2}$ restricted to $\Cbar_{\gamma}$ 
satisfy $\nu_{2}=\frac{[K: \gamma K \gamma^{-1}\cap K]}{[K:K\cap \gamma^{-1}K\gamma]}\,\nu_{1} = \mathrm{mod}(\gamma)\,\nu_{1}$.

\bigskip
\textit{Proof by Mass-Transport Principle}: Denote by $\WW$ the countable discrete set 
$\GGG/K$ and observe that 
$\Xbar$ is isomorphic with the quotient $\GGG\backslash (X\times \WW)$
(where the $\GGG$-action is diagonal). 
$$\left( 
\begin{array}{rcl}
  \GGG\backslash (X\times \WW) & \simeq &K\backslash X= \Xbar\\
(gx,g\gamma K)&\mapsto &K \gamma^{-1} x
\end{array}
\right) 
$$
Denote by $\origin$ the class $K\in 
\GGG/K$ and by $K_{v}$ the stabilizer of $v\in \WW$ in the left 
multiplication $\GGG$-action. In particular, $K_{\origin}=K$.

Two points $\bar{x}_{1}$ and $\bar{x}_{2}$ of 
$\Xbar$ are $\RRbar$-equivalent i{f}{f} they admit 
representatives in $X\times \WW$ with the same first coordinate.
One thus gets an identification of $\RRbar$ with $\GGG\backslash (X\times 
\WW\times \WW)$ (where the $\GGG$-action is diagonal on the three 
coordinates), thanks to the two coordinate-forgetting projections 
(where $\WW_{1}, \WW_{2}$ are two copies of 
$\WW$):
$$
\begin{array}{rrl}
    & & \hskip -10pt\GGG\backslash (X\times \WW_{1})=\Xbar\\
    {}^{pr_{1}}\hskip -10pt & \nearrow  & \\
   \GGG\backslash (X\times \WW_{1}\times \WW_{2}) & & \\
    {}_{pr_{2}}\hskip -10pt & \searrow &  \\
    & & \hskip -10pt\GGG\backslash (X\times \WW_{2})=\Xbar
\end{array}
$$
It becomes equivalent to consider a function $F$ on 
$\RRbar$ or a $\GGG$-invariant function $f$ on 
$X\times \WW\times \WW$. Thus, for 
non-negative functions
$$
\begin{array}{ccl}
    \nu_{1}(F) & = & \displaystyle{\int_{\RRbar} F(\bar{x}_{1}, \bar{x}_{2}) 
    d\nu_{1} = \int_{\Xbar} \sum_{\bar{x}_{2}\sim 
    \bar{x}_{1}}
    F(\bar{x}_{1}, \bar{x}_{2}) d\bar{\mu}(\bar{x}_{1})}
    \\
    & = & \displaystyle{\int_{X}\sum_{v_{2}\in \WW} f(x,\origin,v_{2}) 
    d\mu(x)} \hskip100pt 
\end{array}
$$
while 
$$
\begin{array}{ccl}
    \nu_{2}(F) &= &\displaystyle{\int_{X}\sum_{v_{1}\in \WW} f(x,v_{1},\origin) d\mu(x)}
    \hskip100pt 
\end{array}
$$

On the other hand, the \textit{mass-transport principle} below 
essentially gives the correcting terms for $\nu_{1}$ and $\nu_{2}$ to coincide.
In particular, unimodularity, equivalent to the coincidence of the Haar measures $\Haar(K_{v_{1}})=\Haar(\KK)$ 
for every $v_{1}$, is equivalent to the preservation for $\RRbar$ of the measure 
$\bar{\mu}$.

The \textbf{mass-transport principle}:
$$\displaystyle{\int_{X}\sum_{v_{2}\in \WW}} f(x,\origin,v_{2})\ {\Haar(\KK)}\ d\mu(x)=
\int_{X}\sum_{v_{1}\in \WW} f(x,v_{1},\origin)\ {\Haar(K_{v_{1}})} \ d\mu(x),$$ 
where $\Haar$ is the Haar measure on $\GGG$, is a useful device in invariant percolation theory. 
For details, see \cite{BLPS99a} where I took the 
following two-line proof, credited to W.~Woess.
 Let $\bar{f}(v,v'):=\int_{X} f(x,v,v')\ d\mu(x)$ denote the mean value.
$$
\begin{array}{lll}
    \displaystyle{\sum_{v_{2}\in \WW} \bar{f}(\origin,v_{2})\ \Haar (\KK)}
    & =  \displaystyle{\sum_{v_{2}\in \WW} 
    \bar{f}(\origin,v_{2})\ \Haar(\{g: g\origin=v_{2}\})} 
    & = \displaystyle{\int_{\GGG} \bar{f}(\origin,g\origin)\ d\Haar (g)}
    \\
    \displaystyle{\sum_{v_{1}\in \WW} \bar{f}(v_{1},\origin) \Haar (K_{v_{1}})}
    & =	\displaystyle{\sum_{v_{1}\in \WW}
    \bar{f}(v_{1},\origin) \Haar(\{ \underbrace{g: g v_{1}=\origin}_{\{g: g^{-1}\origin=v_{1}
    \}}\})}
    & = \displaystyle{\int_{\GGG} \bar{f}(g^{-1}\origin,\origin) d\Haar(g) }
\end{array}
$$
And, the last terms are equal, thanks to the $\GGG$-invariance of $f$ 
and $\mu$.
\endofproof

\subsection{Some computation}\label{sect: some computation}
Assume that $\HH$ is transitive and unimodular. 

\noindent
Recall that $\GG_{x}$ denotes, for $x\in X$, the cluster of the 
vertex $\origin$ in the subgraph $\pi(x)$.
Let $p_{x}:C_{(2)}^{1}(\GG)\to\dzero \HD(\GG_{x})$ be
the orthogonal projection from the space of $\ell^2$ cochains of $\GG$ 
to the image, under the coboundary $\dzero$, of $\HD(\GG_{x})$ in $C_{(2)}^{1}(\GG)$.
Denote by $\mathbf{1}_{e_{1}}, \mathbf{1}_{e_{2}}, \cdots, 
\mathbf{1}_{e_{n}}\in C_{(2)}^{1}(\GG)$ the characteristic functions of the
edges $e_{1}, e_{2}, \cdots, e_{n}$ adjacent to the base point 
$\origin$.

\begin{proposition} \label{prop: computation of L2 Betti numbers}
    Let $Y^{\bullet}$ be a non-null $\RRcluster$-saturated
	Borel subset of $X^{\bullet}$, and $Y$ its preimage in $X$.
    The first $L^{2}$ Betti number of the restricted measure equivalence relation
    $\RRcluster_{\vert Y^{\bullet}}$ on $(Y^{\bullet},\mu^{\bullet}_{Y^{\bullet}})$ equals
    $$\beta_{1}(\RRcluster_{\vert 
    Y^{\bullet}},\mu^{\bullet}_{Y^{\bullet}})=
    \frac{1}{2 \mu(Y)}
    \sum_{i=1}^{n} 
    \int_{Y}\scal{p_{y}(\mathbf{1}_{e_{i}})}{\mathbf{1}_{e_{i}}} d\mu(y).$$
\end{proposition}
Here, $\mu^{\bullet}_{Y^{\bullet}}$ is of course the normalized 
restriction of $\mu^{\bullet}$ to $Y^{\bullet}$. The $\frac{1}{2}$ terms just reflects that, the graph 
being transitive, each edge is counted twice: once per endpoint, 
while $\frac{1}{\mu(Y)}$ is just designed to normalize.

\medskip
\proof We use first Theorem~\ref{thm: Betti 1 and graphings} stating 
that $\beta_{1}(\RRcluster_{\vert Y^{\bullet}},\mu_{Y^{\bullet}})=
\dim_{\RRcluster_{\vert Y^{\bullet}}}\int_{Y^{\bullet}} \dzero(\HD(\GG_{y^{\bullet}}))
\ d\mu^{\bullet}_{Y^{\bullet}}(y^{\bullet})$; second the definition of the dimension 
(see \cite[Prop.~3.2 (2)]{Gab02}): A measurable labeling (see Proposition~\ref{prop:quotient is standard})
$e_1^{\bullet}, e_2^{\bullet},
\cdots, e_n^{\bullet}$ of the edges adjacent to $y^{\bullet}$ leads to measurable vector fields 
$y^{\bullet}\mapsto \mathbf{1}_{e_i^{\bullet}}\in C_{(2)}^{1}(\GG_{y^{\bullet}})$ that define a 
family of fields of representative (in the sense of \cite{Gab02}), except that each edge 
is represented twice 
(the additionnal difficulty of 
a possible loop in $\GG$ is dismissed by the fact that it would give a vector orthogonal to $\dzero(\HD(\GG_{y^{\bullet}})$).
Third, we use the relation between the objects with and 
without a $\bullet$ sign. 
\endofproof
\begin{remark}
In case $Y^{\bullet}$ is not $\RRcluster$-saturated, the families of fields of representative are more delicate to describe.
However, Corollary~5.5 of \cite{Gab02}, for induction on Borel subsets (see fact~2, 
Subsection~\ref{subsect: Harmonic Dirichlet Functions and Clusters for Cayley 
Graphs}), applied to $Y^{\bullet}$ and its $\RRcluster$-saturation 
leads to the same formula except that the domain of integration is now the $\HH$-saturation $\HH Y$ of $Y$:
$$\beta_{1}(\RRcluster_{\vert 
    Y^{\bullet}},\mu^{\bullet}_{Y^{\bullet}})=
    \frac{1}{2 \mu(Y)}
    \sum_{i=1}^{n} 
    \int_{\HH Y}\scal{p_{y}(\mathbf{1}_{e_{i}})}{\mathbf{1}_{e_{i}}} d\mu(y).$$
\end{remark}

\medskip
\noindent
Observe that the quantity on the right in the above Proposition~\ref{prop: computation of L2 Betti numbers}
in fact only depends on the image $\Omega_{Y}:=\pi(Y)$ in $\Omega=\{0,1\}^{\EE}$
$$\beta_{1}(\RRcluster_{\vert Y^{\bullet}},\mu^{\bullet}_{Y^{\bullet}})= 
\frac{1}{2 \pi_{*}\mu(\Omega_{Y})}
\sum_{i=1}^{n} 
\int_{\Omega_{Y}}\scal{p_{\omega}(\mathbf{1}_{e_{i}})}{\mathbf{1}_{e_{i}}} 
d\pi_{*}\mu(\omega).$$
Concerning the full equivalence relation, one gets, 
$$
    \beta_{1}(\RRfull,\mu^{\bullet})=\frac{1}{2} 
    \sum_{i=1}^{n} 
    \int_{X}\scal{p_{x}(\mathbf{1}_{e_{i}})}{\mathbf{1}_{e_{i}}} 
    d\mu(x).
$$
But in the case of $\RRfull$, the projection $p_{x}$ doesn't depend 
on $x$: it is just the projection $$p:C_{(2)}^{1}(\GG)\to\dzero 
\HD(\GG)\simeq \HD(\GG)/\Cmath.$$
It follows that 

\begin{proposition}
For the full equivalence relation on 
$X^{\bullet}$:
    \begin{eqnarray}
	\beta_{1}(\RRfull,\mu^{\bullet})=\frac{1}{2} 
	\sum_{i=1}^{n} 
	\scal{p(\mathbf{1}_{e_{i}})}{\mathbf{1}_{e_{i}}}.
    \end{eqnarray}
\end{proposition}
Observe that this quantity doesn't even depend on what happens on $X$ 
nor on the choice of $\HH$, once $\HH$ is unimodular and transitive on 
the vertices: It is an invariant of the graph. 
\begin{definition}\label{def: beta1 of the graph}
Call it the {\bf first 
$\ell^2$ Betti number} of $\GG$ and denote it 
\begin{eqnarray}\label{eq: def of beta(graph)}
\beta_{1}(\GG):=\frac{1}{2} 
    \sum_{i=1}^{n} 
    \scal{p(\mathbf{1}_{e_{i}})}{\mathbf{1}_{e_{i}}}.
\end{eqnarray}
\end{definition}
It is clear\footnote{\label{footn: beta1 and OHD}$
    \scal{p(\mathbf{1}_{e_{i}})}{\mathbf{1}_{e_{i}}}
    = 
    \scal{p^2(\mathbf{1}_{e_{i}})}{\mathbf{1}_{e_{i}}}
    =
    \scal{p(\mathbf{1}_{e_{i}})}{p(\mathbf{1}_{e_{i}})}=0$ i{f}{f}
    $h.p(\mathbf{1}_{e_{i}})=p(h.\mathbf{1}_{e_{i}})=0$
    for all $h\in\HH$ and ${i=1,\cdots,n}$.
}
that $\beta_1(\GG)=0$ if and only if $\GG$ belongs to $\OOO$.

\section{Harmonic Dirichlet Functions and Clusters for Transitive
Graphs} \label{sect: Harmonic Dirichlet Functions and Clusters for Transitive
Graphs}

In this section, we give the proof of Theorem~\ref{theo intro: when are selectable 
clusters in OHD} of the introduction by putting/proving it in  the more 
general context of selectability.

Recall that we have at hand: (1) A locally finite graph $\GG=(\VV,\EE)$ with a  
transitive (on $\VV$) action of a closed group $\HH$ of automorphisms;
(2) An $\HH$-invariant probability measure $\mathbf{P}$ on the set $\Omega=\{0,1\}^{\EE}$
of colorings of $\GG$.
	
\subsection{The action made free}
The $\HH$-action on $\Omega$ being not necessarily free, let's
consider a diagonal $\HH$-action on $X=\Omega\times Z$,
where $Z$ is a standard Borel space with an essentially
free
probability measure-preserving action of $\HH$.  An example of such a
$Z$ is furnished by the Lemma~\ref{lem: there exists a free action of
G} below.

The diagonal action preserves the product measure $\mu$ and is
(essentially) free. The obvious projection $\pi:  
X=\Omega\times Z  \to \Omega$
sends $\mu$ to $\mathbf{P}$ and is $\HH$-equivariant, so that we are 
in the context of Section~\ref{sect:Percolation on transitive graphs}.

\begin{remark}\label{Rem: main point of the proof}
    We could probably avoid the detour by the freeness of the action 
    by defining $L^2$ Betti numbers for groupoids instead of 
    just for equivalence relations, as suggested in 
    \cite[p.103]{Gab02}. Notice that such a study of $L^2$ Betti 
    numbers for measured groupoids has been 
    carried out by R.~Sauer (see \cite{Sau03}), 
    using L\"uck's approach of $\ell^2$ theory.
\end{remark}

\begin{lemma}\label{lem: there exists a free action of G}
    If $\HH$ acts continuously\,\footnote{The stabilizer of a point of 
    $V$ is a closed open subgroup of $\HH$.} faithfully on a discrete countable set 
    $V$, then the diagonal action of $\HH$  on 
    ${\check{\Omega}}=\left(\{0,1\}^{V}\right)^{\Nmath}$ is continuous, preserves the 
    Bernoulli measure (product of equiprobabilities on $\{0,1\}$), 
    and is essentially free.
\end{lemma}
\proof Enumerate the elements of $V$: $v_{1},v_{2},\ldots, v_{n},\ldots$ and 
denote by  ${\check{\Omega}}_{i,j}$ the subset of points of ${\check{\Omega}}$ that 
are fixed by an element of $\HH$ which sends $v_{i}$ to $v_{j}$. 
This subset satisfying infinitely many equations: $\omega(v_{i},l)=\omega(v_{j},l)$
for each coordinate $l\in \Nmath$, has thus measure 0. The set of points $\omega$ with 
a non-trivial stabilizer is contained in the countable union of the 
${\check{\Omega}}_{i,j}$; it has measure zero.
\endofproof

\subsection{Selectability}\label{sect:about 
selectability}

Let $\GG$ be a locally finite transitive graph and $\HH$
a closed transitive subgroup  of $\mathrm{Aut}(\GG)$.
Equip $\{0,1\}^{\VV}$ with the natural action of $\HH$ induced
by its action on $\VV$.
\begin{definition}\label{def: Delta select. cluster} 
Let $\mathbf{P}$ be an $\HH$-invariant percolation on $\GG$.
Let $\Delta$ be a closed subgroup of $\HH$.
\begin{description}
\item
A $\Delta$-equivariant \textbf{selected cluster} on a $\Delta$-invariant $\mathbf{P}$-non-null Borel subset $D\subset\{0,1\}^{\EE}$ is a 
$\Delta$-equivariant measurable map $c:D\to\{0,1\}^{\VV}$, 
such that $c(\omega)$ is the (characteristic function of
the) vertex set of one cluster $C(\omega)$ of $\omega$.
\item
A $\Delta$-equivariant \textbf{virtually selected cluster} 
on a $\Delta$-invariant $\mathbf{P}$-non-null Borel subset $D\subset\{0,1\}^{\EE}$ is a 
$\Delta$-equivariant measurable map $c:D\to\{0,1\}^{\VV}$
such that $c(\omega)$ is the (characteristic function of
the) vertex set of (the union of) finitely many clusters $C_{1}(\omega), 
C_{2}(\omega), \cdots, C_{n(\omega)}(\omega)$ of $\omega$.
\end{description}
\end{definition}

\begin{example}
If almost every subgraph $\omega$ has a unique 
infinite cluster, 
assigning to $\omega$ this infinite cluster defines an $\HH$-equivariant 
selected cluster.
And similarly, if almost every subgraph $\omega$ has  finitely many
infinite clusters, 
assigning to $\omega$ these infinite clusters defines an 
$\HH$-equivariant virtually selected cluster.
\end{example}

Let $(X,\mu)$ be a standard Borel probability space together with 
\begin{itemize}
    \item a probability measure-preserving (p.m.p.) action of 
    $\HH$, which is free, and
    \item an $\HH$-equivariant Borel map (field of graphs) $\pi:X\to \{0,1\}^{\EE}$
\end{itemize}
so that  our situation fits in the general context of
Section~\ref{sect:Percolation on transitive graphs}.

\begin{proposition}\label{prop: Rfull=Rcluster on selected set}
  The following are equivalent:\\
(1) The invariant percolation $\pi_*(\mu)$ admits a 
 $\HH$-equivariant selected cluster,\\
(2) There is a non-null Borel subset $\TT$ of
    $X^{\bullet}$ to which the restrictions of $\RRfull$ and $\RRcluster$
    coincide: $\RRfull_{\vert \TT}=\RRcluster_{\vert 
    \TT}$.

Moreover, $\TT$ can be taken to be the image in $X^{\bullet}$ of the set of those $x\in X$
whose selected cluster contains the base point $\origin$.
\end{proposition}

\proof  
The selected cluster
$c:\{0,1\}^{\EE}\to \{0,1\}^{\VV}$ selects by composition by $\pi$
one connected component of the splitting of each $\Lfull$-leaf into its
$\Lcluster$-components. The union of these selected leaves intersects
the transversal $X^{\bullet}$ along a Borel subset $\TT$ which is 
characterized as the image in $X^{\bullet}$ of the set of $x\in X$ such that
the selected cluster $c(\pi(x))$ contains the base point $\origin$. 
Two
$\RRfull$-equivalent points in $\TT$ belong to the same $\Lfull$-leaf
and both belong to THE selected
$\Lcluster$-leaf; they are thus $\RRcluster$-equivalent.

  Conversely, if
$\TT$ is a Borel subset of $X^{\bullet}$ to which the restrictions of
$\RRfull$ and $\RRcluster$ coincide, then it selects an 
$\Lcluster$-leaf in each $\Lfull$-leaf meeting $\TT$.
One can assume that $\TT$ is $\RRcluster$-saturated (two points that are 
$\RRfull$-equivalent and $\RRcluster$-equivalent to some point in $\TT$
have to be $\RRcluster$-equivalent).
Let $\check{T}$ be the preimage of $\TT$ in $X\times \VV$.
It is an $\HH$-invariant subset, whose projection $T$ in $X$
is non-null and $\HH$-invariant,
and whose intersection with each fiber $(\{x\}\times\VV)\cap\check{T}$
is a cluster $\check{c}(x)$ of $\pi(x)$. This defines an $\HH$-equivariant map 
$\check{c}:T\to \{0,1\}^{\VV}$. 
Now the set $(X\times \{\origin\})\cap\check{T}$ 
once projected in $X$ corresponds to those
$x\in T$ for which $\check{c}(x)$ is the cluster of the base point $\origin$.
This shows that the cluster $\check{c}(x)$ only depends on $\pi(x)$.
Moreover,  $\pi(T)$ is non-null for the measure  $\pi_*\mu$, 
so that the map $\check{c}$ induces an $\HH$-equivariant 
selected cluster on $\pi(T)\subset\{0,1\}^{\EE}$. 
\endofproof

\medskip
Recall that a subrelation ${\mathcal S}$ has \emph{finite index} in
$\RR$ if each $\RR$-class splits into finitely many ${\mathcal
S}$-classes.
The same kind of argument as above shows:
\begin{proposition}\label{prop: [Rfull:Rcluster]
finite on virt.  selected set}
  The following are equivalent:\\
(1) The invariant percolation $\pi_*(\mu)$ admits a 
 $\HH$-equivariant virtually selected cluster,\\
(2) There is a non-null Borel subset $\TT$ of
    $X^{\bullet}$ to which the restriction of $\RRcluster$ has finite 
    index in the restriction of $\RRfull$: $[\RRfull_{\vert \TT}:\RRcluster_{\vert 
    \TT}]<\infty$.

Moreover, $\TT$ can be taken to be the image in $X^{\bullet}$ of the set of those $x\in X$
for which one of the selected clusters contains the base point $\origin$.
\end{proposition}

Observe that if $\GG$ is not a finite graph, then $\RRfull$ has 
infinite classes and the $\HH$-equivariant virtually selected clusters are (almost) all infinite.

\begin{remark}
The main result of \cite{LS99} (indistinguishability of the infinite clusters) implies that Bernoulli percolation
in the nonuniqueness phase admits no $\HH$-equivariant virtually selected clusters.
\end{remark}

\medskip
\begin{remark}
{\em
Let $\Delta$ be a closed subgroup of $\HH$ that contains 
the stabilizer $\KK$ of $\origin$. One can define a notion of 
full $\Delta$-equivalence relation $\RRfull_{\Delta}\subset \RRfull$
and cluster $\Delta$-equivalence relation: the intersection
$\RRcluster_{\Delta}=\RRcluster\cap \RRfull_{\Delta}$.

\medskip
\begin{tabular}{lp{14cm}}
    $\RRfull_{\Delta}$: &
    \emph{Two points $x^{\bullet},y^{\bullet}$ are 
	$\RRfull_{\Delta}$-equivalent if 
    and only if 
    they admit two representatives $x,y$ in $X$ for which there exists 
    $\delta\in \Delta$ such that $\delta x=y$, i{f}{f} 
    any of their representatives are in the same $\Delta$-orbit.}
\end{tabular}
\medskip

\noindent
Similarly, for the cluster $\Delta$-equivalence relation:

\medskip
\begin{tabular}{lp{14cm}}
    $\RRcluster_{\Delta}$: &
    \emph{Two points $x^{\bullet},y^{\bullet}$ are $\RRcluster_{\Delta}$-equivalent if and only if 
    they admit two representatives $x,y$ in $X$ for which there exists 
    $\delta\in \Delta$ such that $\delta x=y$  and the vertices 
$\origin, \delta^{-1} \origin$ 
    are in the same cluster of $\pi(x)$.}
\end{tabular}
\medskip

In case $\GG$ is the Cayley graph of a discrete group $\GGamma$
(i.e.\  $\KK=\{1\}$ and $X^{\bullet}=X$) then $\RRfull_{\Delta}$
is just the equivalence relation defined by the $\Delta$-action
on $X$, while $\RRcluster_{\Delta}=\RRcluster\cap \RRfull_{\Delta}$ is 
just defined by: 
$(x,y)\in \RRcluster_{\Delta}$ i{f}{f} there exists 
$\delta\in \Delta$ such that $\delta x=y$  and the vertices 
$\origin, \delta^{-1} \origin$ are in the same cluster of $\pi(x)$.

Exactly along the same arguments as above, one can show that the following are equivalent:\\
{\em 
(1) The invariant percolation $\pi_*(\mu)$ admits a 
 $\Delta$-equivariant selected cluster,\\
(2) There is a non-null Borel subset $\TT$ of
    $X^{\bullet}$ to which the restrictions of $\RR_{\Delta}$ and $\RRcluster_{\Delta}$ 
	coincide: $\RRfull_{\Delta\vert T}=\RRcluster_{\Delta\vert T}$. \\
And similarly, with 	
finite index, for the virtual notion.
}

The lamination interpretation of these equivalence relations goes as follows:\\
Consider first the space $X\times \GG$ and divide out by $\Delta$
to get the laminated space $\Lfull_{\Delta}$.
Consider now the transversal $\Delta\backslash
(X\times \Delta\origin)\subset \Delta\backslash
(X\times \VV)$, which is naturally isomorphic 
with $X^{\bullet}=\KK\backslash X$ since $\Delta$ contains $\KK$,
and the equivalence relation defined on it 
by ``belonging to the same $\Lfull_{\Delta}$-leaf''. This is the full $\Delta$-equivalence relation $\RRfull_{\Delta}\subset \RRfull$ and it appears as 
the image in $X^{\bullet}=\HH\backslash (X\times \VV)$ of the equivalence relation defined 
by the $\Delta$-action on $X$.
Just like in Section~\ref{sect:Percolation on Cayley Graphs}, use now $\pi$  to get a coloring on the leaves.
Define $\Lcluster_{\Delta}$ as the sub-laminated space where
the $0$-colored edges are removed and $\RRcluster_{\Delta}$ 
as the subrelation of $\RRfull_{\Delta}$ induced on $\Delta\backslash
(X\times \Delta\origin)$ by ``belonging to the same 
$\Lcluster_{\Delta}$-leaf''.
}
\end{remark}

\subsection{Selected Clusters and Harmonic Dirichlet Functions}
\label{selected Clusters and Harmonic Dirichlet Functions}

The connections between selected clusters and harmonic Dirichlet 
functions is very simple:
\begin{theorem}\label{th: virt select clusters in/out OHD}
    Assume $\pi_*\mu$ admits an $\HH$-equivariant 
virtually selected cluster.  Assume that
    the closed subgroup $\HH$ is unimodular.
    If $\GG$ belongs to $\OOO$, then $\mu$-a.e.\ virtually selected
    cluster belongs to $\OOO$.
    If $\GG$ doesn't belong to $\OOO$, then $\mu$-a.e.\ virtually
    selected cluster doesn't belong to $\OOO$.
\end{theorem}

\proof
Thanks to unimodularity, the associated equivalence relations 
$\RRfull$ and $\RRcluster$ are measure-preserving (Th.~\ref{th: measure-preserved iff unimodular}).

Start with the case of a selected cluster.
Denote by $c$ the $\HH$-equivariant selected cluster $c:\{0,1\}^{\EE}\to \{0,1\}^{\VV}$
and $C=c\circ \pi:X\to \{0,1\}^{\VV}$.
Let $T$ be the Borel subset of $x \in X$ where the selected 
cluster $C(x)$ contains $\origin$ and let $\TT$ be its image in $X^{\bullet}$.

The following are equivalent:
\begin{enumerate}
    \item  $\GG$ is in $\OOO$ \label{enum: unimod-G dans OHD}

    \item  $\beta_{1}(\RRfull,\mu^{\bullet})=0$ \label{enum: unimod-beta1 de Sfull 
    vanishes}
    
    \item  $\beta_{1}(\RRfull_{\vert \TT},\mu^{\bullet}_{\TT})=0$ \label{enum: unimod-beta1 de Sfull IY vanishes}
    
    \item  $\beta_{1}(\RRcluster_{\vert \TT},\mu^{\bullet}_{\TT})=0$ \label{enum: unimod-beta1 
    de Scluster vanishes}

    \item  for $\mu^{\bullet}$-almost every $x^{\bullet}\in \TT$, the graph 
    $\Lcluster_{x^{\bullet}}$ is in $\OOO$
    \label{enum: unimod-a.s. Sigmacluster in OHD}
    
    \item  for $\mu$-almost every $x\in T$, the selected 
    cluster $C(x)$ is in $\OOO$ \label{enum: unimod-a.s. upstair clusters in OHD}
\end{enumerate}

Theorem~\ref{thmintro: Betti 1 and graphings}
applied to the field of graphs $x^{\bullet}\mapsto 
\Lfull_{x^{\bullet}}\simeq \GG$ 
(example~\ref{examp: fields of graphs for transitive}) gives the
equivalence \ref{enum: unimod-G dans OHD}$\iff$\ref{enum: unimod-beta1 de Sfull 
vanishes}. When applied to the restriction of that field to $\TT$
(example~\ref{examp: smooth fields of graphs+restrictions}), it gives 
\ref{enum: unimod-G dans OHD}$\iff$\ref{enum: unimod-beta1 de Sfull IY vanishes}.
Observe that the equivalence \ref{enum: unimod-beta1 de Sfull
vanishes}$\iff$\ref{enum: unimod-beta1 de Sfull IY vanishes} is also
an application of \cite[Cor.~5.5]{Gab02}.
Since $\beta_{1}$ is an invariant of the equivalence relation, 
\ref{enum: unimod-beta1 de Sfull IY vanishes}$\iff$\ref{enum: unimod-beta1 
de Scluster vanishes} is deduced from the coincidence 
$\RRcluster_{\vert \TT}=\RRfull_{\vert \TT}$
(prop.~\ref{prop: Rfull=Rcluster on selected set}).
Theorem~\ref{thmintro: Betti 1 and graphings}, applied to the field of graphs $x^{\bullet}\mapsto 
\Lcluster_{x^{\bullet}}$ 
(example~\ref{examp: fields of graphs for transitive}) restricted to 
$\TT$, shows the equivalence \ref{enum: unimod-beta1 
de Scluster vanishes}$\iff$\ref{enum: unimod-a.s. Sigmacluster in OHD}.
Each $\Lcluster_{x^{\bullet}}$ being isomorphic to the cluster 
$\pi(x)(\origin)=C(x)$ of any of its representatives, and $\mu^{\bullet}$ being 
the push-forward of $\mu$, one deduces the last equivalence.

Let $T'\subset T$ be the Borel subset of points such that the
selected cluster is in $\OOO$.  If $T'$ is non-null, then the
above arguments applied to $T'$ shows that $\GG$ belongs to $\OOO$ and
thus $T'=T$ a.s.\ This implies that in case $\GG$ does not belong to
$\OOO$, then for $\mu$-almost every $x\in T$, the selected cluster
$C(x)$ is not in $\OOO$.  
\bigskip

For the case of a virtually selected cluster, partition first $T$
as $\coprod T_{n}$ according to how many clusters are selected.  On
$\TT_{n}$ the $\RRfull_{\vert \TT}$-classes decompose into $n$
$\RRcluster_{\vert \TT}$-classes (prop.~\ref{prop: [Rfull:Rcluster]
finite on virt.  selected set}).  Then just replace the argument in the proof
of the equivalence \ref{enum: unimod-beta1 de Sfull IY
vanishes} $\iff$\ref{enum: unimod-beta1 de Scluster vanishes} above
by Proposition~5.11 of \cite{Gab02}, asserting that
$\beta_{1}(\RRcluster_{\vert \TT})=n \beta_{1}(\RRfull_{\vert \TT})$.
\endofproof

\section{Nonuniqueness Phase and Harmonic Dirichlet Functions}
\label{sect: Non uniq. pahse Harm. Dir. Funct}
This section is concerned with a comparison between  two invariant percolations.
Its main goal is to prove Theorem~\ref{th: general pc<pu+ couplings}, which implies 
both Theorem~\ref{thintro: pc<pu for G not in OHD} 
(Corollary~\ref{cor: pc<pu for G not in OHD}) and 
Theorem~\ref{thintro: beta-1 and pc,pu for Random-Cluster} 
(Corollary~\ref{cor: `pc<pu' for G not in OHD and RC model}) of the introduction.

Consider a unimodular transitive group $\HH$ of automorphisms of $\GG$ and
two $\HH$-invariant percolations $\mu_{1}$ and $\mu_{2}$ on $\GG$.
Recall that given two $\HH$-invariant percolations on $\GG$, an 
\textbf{$\HH$-equivariant coupling} is a p.m.p.\ $\HH$-action on
a standard probability measure space $(X,\mu)$ with two $\HH$-equivariant maps
$$
\begin{array}{rcl}
    & (X,\mu) &\\
    {}^{\pi_{1}}\swarrow & & \searrow^{\pi_{2}}\\
    (\{0,1\}^{\EE},\mu_{1}) & & (\{0,1\}^{\EE}, \mu_{2})
\end{array}
$$ pushing $\mu$ to $\mu_{i}$ respectively, i.e.\   $\pi_{1*}\mu=\mu_{1}$ 
and $\pi_{2*}\mu=\mu_{2}$.

To $\pi_{1}$ and $\pi_{2}$ correspond two laminations 
$\Lcluster_{1}$ and $\Lcluster_{2}$ (both sub-laminations of $\Lfull$) with transversal 
$X^{\bullet}=\KK\backslash X$ and two cluster equivalence 
relations $\RRcluster_{1}$ and $\RRcluster_{2}$ (both subrelations of $\RRfull$).
\begin{remark}
An invariant coupling always exists since the product space with the 
product measure and the diagonal action will do.
In case the coupling witnesses
 a stochastic domination\,\footnote{i.e.\ 
$\mu\{x: \pi_{1}(x)\leq \pi_{2}(x)\}=1$, i.e.\  
for $\mu$ a.e.\ $x\in X$ and for every edge $e\in \EE$, if 
$\pi_{1}(x)(e)=1$, then $\pi_{2}(x)(e)=1$.}
(see for instance \cite[sect. 2.2]{HJL02a}), 
then $\RRcluster_{1}$ fits into $\RRcluster_{2}$: $\RRcluster_{1}\subset \RRcluster_{2}$.
\end{remark}

\begin{theorem}\label{th: general pc<pu+ couplings}
 Let $\GG$ be a unimodular transitive locally finite graph.
 Let $\HH$ be a unimodular transitive group of automorphisms of 
    $\GG$, let $(X,\mu)$ be an $\HH$-equivariant coupling
    between two $\HH$-invariant percolations $\mu_{1}, \mu_{2}$.
	Assume that\\
	\indent 1. $\mu_{1}$-a.e. cluster belongs to $\OOO$,\\
	\indent 2. $\pi_2$ has an $\HH$-equivariant selected cluster defined on a 
 non-null set\,\footnote{For instance, if
 $\mu_{2}$ has a non-null set of subgraphs with exactly one 
   infinite cluster.},\\
then
	$$\displaystyle{\beta_{1}(\GG)\leq \frac{1}{2} \hskip-20pt
	\sum_{{\begin{array}{c}
    \mbox{\footnotesize{\em edges\ }} e\\
    \mbox{\footnotesize{\em adjacent\ to\ }} \origin
    \end{array}}}\hskip-20pt
	\mu\Bigl(\pi_{2}(e)=1 \mathrm{\ and\ } \pi_{1}(e) =0 \Bigr)
	}.$$
\end{theorem}
Here, ${\beta_{1}(\GG)}$ is the invariant of the graph introduced in Definition~\ref{def: beta1 of the graph}.
It is strictly positive if and only if $\GG$ doesn't belong to $\OOO$.

The main ingredient in the proof of the theorem will be the following 
useful result. Here, graphing may be understood as oriented or 
unoriented (see Section~\ref{sect: background Measured Equivalence Relations}).
\begin{theorem}\label{th: R2 = R1 vee psi2 + beta1 and cost}
Let $\RR_1$ be a p.m.p.\ equivalence relation on the standard Borel space $(X,\mu)$. Let $\Psi_{2}$ be a 
p.m.p.\ graphing
and let $\RR_2=\RR_1\vee \Psi_{2}$ be the equivalence relation generated by $\RR_1$ and $\Psi_{2}$.
Then $$\beta_1(\RR_2)-\beta_0(\RR_2)\leq \beta_1(\RR_1)-\beta_0(\RR_1)+\mathrm{cost}(\Psi_{2}).$$
\end{theorem}
\begin{remark}
    There is no continuity in the other direction: 
    Think in $\RR_i$ given by a free action of $\Gamma_1=\FF_2$ and $\Gamma_2=\FF_2\times \Zmath$.
    Then $\RR_2$ can be generated from $\RR_1$ by adding a graphing $\Psi_2$ of arbitrarily small cost.
However, $\beta_1(\RR_2)-\beta_0(\RR_2)=0$ while 
    $\beta_1(\RR_1)-\beta_0(\RR_1)=1$.
\end{remark}

\proof
The proof is just an adaptation of the proof of the Morse inequalities 
(see [Gab02, sect. 4.4, p.137]).
Let $\bar{\Sigma}_1$ be a \emph{simply connected} smooth simplicial $\RR_1$-complex, with a big enough $0$-skeleton: $\RR_1\subset \bar{\Sigma}_1$. Then
$(\beta_1-\beta_0)(\RR_1)=(\beta_1-\beta_0)(\bar{\Sigma}_1)$. Let $(\bar{\Sigma}{}^{n}_1)_n$ be an increasing
sequence of ULB smooth simplicial $\RR_1$-complexes that exhausts $\bar{\Sigma}_1$, with say $\RR_1\subset\bar{\Sigma}_1^{0}$.\\ 
Let $\bar{\bar{\Sigma}}_1$ and
$\bar{\bar{\Sigma}}{}^{n}_1$ be the corresponding smooth
$\RR_2$-complexes ([Gab02, sect. 5.2, p.140]): $\RR_2\subset \bar{\bar{\Sigma}}_1^{0}$ and moreover, the reciprocity formula gives
$$(\beta_1-\beta_0)(\RR_1)= \lim_n (\beta_1-\beta_0)(\bar{\Sigma}{}^{n}_1, \RR_1) = \lim_n
(\beta_1-\beta_0)(\bar{\bar{\Sigma}}{}^{n}_1, \RR_2).$$
Define $\bar{\bar{\Sigma}}_2$ and $\bar{\bar{\Sigma}}{}^{n}_2$ by just adding to 
$\bar{\bar{\Sigma}}_1$ and $\bar{\bar{\Sigma}}{}^{n}_1$ the $\RR_2$-field of graphs associated to $\Psi_{2}$ on $\RR_2$.\\
Claim: $(\beta_1-\beta_0)(\bar{\bar{\Sigma}}{}^{n}_2)\leq (\beta_1-\beta_0) (\bar{\bar{\Sigma}}{}^{n}_1)+ \mathrm{cost}(\Psi_{2})$, $\forall n$.\\
This is an immediate computation by mimicking the first 4 lines of [Gab02, sect. 4.4, p.137] (read $b_0=0$ instead of $b_{-1}=0$, there).\\
Now, by letting $n$ tend to $\infty$, $(\beta_1-\beta_0)(\bar{\bar{\Sigma}}_2)\leq (\beta_1-\beta_0) (\bar{\bar{\Sigma}}_1)+ \mathrm{cost}(\Psi_{2})=(\beta_1-\beta_0)
(\RR_1)+ \mathrm{cost}(\Psi_{2})$.\\ 
By definition, $\bar{\bar{\Sigma}}_2$ is connected, so that $\beta_0(\bar{\bar{\Sigma}}_2)=\beta_0(\RR_2)$ [Gab02,
3.14].
On the other hand, $\beta_1(\RR_2)\leq \beta_1(\bar{\bar{\Sigma}}_2)$ [Gab02, 3.13 and 3.14].
This proves Th.~\ref{th: R2 = R1 vee psi2 + beta1 and cost}.
\endofproof

\bigskip
{\textsc{Proof}} (of Theorem~\ref{th: general pc<pu+ couplings}):
Denote by $\check{\RR}=\RRcluster_{1}\vee \RRcluster_{2}$ the equivalence relation 
generated by $\RRcluster_{1}$ and $\RRcluster_{2}$. It is also the relation defined 
on $X^{\bullet}$ by the union lamination 
$\Lcluster_{1}\cup\Lcluster_{2}$.\\
Let $Y\subset X^{\bullet}$ be a $\mu^{\bullet}$-non-null 
$\check{\RR}$-saturated Borel subset.
It is clear from the lamination description that the restriction 
$\check{\RR}_{\vert Y}$ is generated by the restriction $\RRcluster_{1\vert 
Y}$ together with the graphing $\Psi_{2}$ consisting of edges of 
$\Lcluster_{2}\setminus \Lcluster_{1}$ with both endpoints in $Y$ (one 
endpoint in $Y$ implies the other one in $Y$, by saturation):
$$\check{\RR}_{\vert Y}=\RRcluster_{1\vert Y}\vee \Psi_{2}.$$
The  cost of $\Psi_{2}$ (again, $\frac{1}{2}$ just reflects that each edge is counted twice, 
while $\frac{1}{\mu^{\bullet}(Y)}$ is just designed to normalize) is bounded above by 
$$\mathrm{cost}(\Psi_{2})\leq \frac{1}{2 \mu^{\bullet}(Y)}\hskip-20pt
\sum_{{\begin{array}{c}
    \mbox{\footnotesize{edges\ }} e\\
    \mbox{\footnotesize{adjacent\ to\ }} \origin
    \end{array}}}
\hskip-20pt
\mu\Bigl(\pi_{2}(e)=1 \mathrm{\ and\ } \pi_{1}(e) =0 \Bigr).$$
The above Theorem~\ref{th: R2 = R1 vee psi2 + beta1 and cost} 
gives:
$$\beta_{1}(\check{\RR}_{\vert Y})-\beta_{0}(\check{\RR}_{\vert 
Y})\leq
\beta_{1}(\RRcluster_{1\vert Y})-\beta_{0}(\RRcluster_{1\vert Y})+\frac{1}{2 \mu^{\bullet}(Y)} 
\hskip-20pt
\sum_{{\begin{array}{c}
    \mbox{\footnotesize{edges\ }} e\\
    \mbox{\footnotesize{adjacent\ to\ }} \origin
    \end{array}}}\hskip-20pt
\mu\Bigl(\pi_{2}(e)=1 \mathrm{\ and\ } \pi_{1}(e) =0 \Bigr).$$
If $Z$ is a measurable non-null subset where $\origin$ belongs to the 
$\HH$-equivariant selected cluster (for instance, the unique infinite cluster),
then $\RRcluster_{2\vert Z}=\RRfull_{\vert Z}$. So that its 
$\check{\RR}$-saturation $Y$ satisfies $\check{\RR}_{Y}=\RRfull_{Y}$ and thus
(see fact~2, 
Subsection~\ref{subsect: Harmonic Dirichlet Functions and Clusters for Cayley 
Graphs})
$$\mu^{\bullet}(Y)\beta_{1}(\check{\RR}_{\vert Y})=\mu^{\bullet}(Y)
\beta_{1}(\RRfull_{\vert Y})=\beta_{1}(\RRfull),$$
which coincides by definition with the quantity $\beta_{1}(\GG)$ introduced in
Section~\ref{sect: some computation}, Definition~\ref{def: beta1 of the graph}.\\
Recall that $\beta_{1}(\GG)=0$ i{f}{f} $\GG\in \OOO$. The graph $\GG$ being infinite, 
$\beta_{0}(\check{\RR}_{\vert 
Y})=\beta_{0}(\RRfull_{\vert Y})=0$, so that 
$$0<\beta_{1}(\GG)\leq
\mu^{\bullet}(Y) \bigl(\beta_{1}(\RRcluster_{1\vert Y})-\beta_{0}(\RRcluster_{1\vert Y})\bigr)+\frac{1}{2} 
\hskip-20pt
\sum_{{\begin{array}{c}
    \mbox{\footnotesize{edges\ }} e\\
    \mbox{\footnotesize{adjacent\ to\ }} \origin
    \end{array}}}\hskip-20pt
\mu\Bigl(\pi_{2}(e)=1 \mathrm{\ and\ } \pi_{1}(e) =0 \Bigr).$$
Now, the assumption (1) of Theorem~\ref{th: general pc<pu+ couplings} 
is designed (since $\beta_{1}(\RRcluster_{1\vert Y})=0$ by Th.~\ref{thmintro: Betti 1 and graphings})
 to ensure that
$$\beta_{1}(\RRcluster_{1\vert Y})-\beta_{0}(\RRcluster_{1\vert Y})\leq 
0$$ and this finishes its proof.
\endofproof
\subsection{Application to Bernoulli Percolation}
\label{subsect: Application to Bernoulli Percolation}
\begin{corollary}{(Th.~\ref{thintro: pc<pu for G not in OHD})}
\label{cor: pc<pu for G not in OHD}    
	Let $\GG$ be a unimodular transitive locally finite graph.
    If $\GG$ doesn't belong to $\OOO$, then the nonuniqueness phase 
    interval of Bernoulli percolation has non-empty interior:
    $$p_{c}(\GG) < p_{u}(\GG)$$
	More precisely, 
	$$0<\beta_{1}(\GG)\leq \frac{1}{2} (\mathrm{degree\ of\ } \GG) (p_{u}(\GG)-p_{c}(\GG)).$$
\end{corollary}

{\textsc{Proof}}:
The standard coupling 
$([0,1]^{\EE}, \otimes\mathrm{Lebesgue})\overset{\pi_p}{\rightarrow} (\{0,1\}^{\EE}, \mathbf{P}_p)$ 
(see for example Section~\ref{examp: standard coupling}) provides a family of countable equivalence relations $\RRcluster_{p}$ (on the quotient space 
$\KK\backslash [0,1]^{\EE}$, once given a closed unimodular transitive group of automorphisms of $\GG$).
For $s<p_{c}$, the equivalence classes of the 
cluster equivalence relation are a.s.\ finite, thus $\mu_{s}$-a.e.\ 
cluster belongs to $\OOO$.
The right-hand quantity of 
Theorem~\ref{th: general pc<pu+ couplings}
with $\mu_{1}=\mu_{s}$ and $\mu_{2}=\mu_{t}$, for $t$ in the uniqueness phase,
is 
$$\displaystyle{\frac{1}{2} \hskip-20pt
\sum_{{\begin{array}{c}
    \mbox{\footnotesize{edges\ }} e\\
    \mbox{\footnotesize{adjacent\ to\ }} \origin
    \end{array}}}\hskip-20pt
\mu\Bigl(\pi_{2}(e)=1 \mathrm{\ and\ } \pi_{1}(e) =0 \Bigr)  \ =\ 
\frac{1}{2} (\mathrm{degree\ of\ } \GG) (t-s)}.$$
One concludes by continuity, by letting $s$ tend to $p_c(\GG)$ and 
$t$ tend to $p_u(\GG)$.

Observe that one could have applied Theorem~\ref{th: general pc<pu+ couplings}
directly with $s=p_c(\GG)$: there is almost surely no infinite cluster
at $p_c(\GG)$ \cite[Th.~1.3]{BLPS99a}.
\endofproof

\begin{remark}
\em{While the above corollary extends to 
unimodular quasi-transitive locally finite graphs, 
it is unknown whether the unimodularity assumption may be removed.
On the other hand, the removal of any transitivity assumption makes it false
since R.~Lyons and Y.~Peres showed (personal communication: I want to
thank them allowing me to reproduce their description here) that the following
graph $\GG$ doesn't belong to $\OOO$ but on the other hand, 
the set of parameters $p$ in Bernoulli percolation 
doesn't admit any interval of nonuniqueness.

Denote by $\GG_m$ the graph obtained from the lattice $\Zmath^2$ by replacing
each edge of $\Zmath^2$ by $m$ paths of length 2. We fix $m$ large enough
so that $p_c(\GG_m)<p_c(\Zmath^3)$. Denote by $\GG_m(k)$ a $k$-by-$k$ square in
$\Zmath^2$, with each edge replaced by $m$ paths of length 2.
Now consider two copies of $\Zmath^3$ (call them $\GG'$ and $\GG''$)
that we will join in countably many corresponding places $(x_i, y_i)\in \GG'\times \GG''$,
with density $0$ in both $\GG'$ and $\GG''$,
using graphs $\GG_m(k_i)$ as follows: position $\GG_m(k_i)$
with one corner at $x_i$ and another corner at $y_i$.
We make $k_i$ grow fast enough so that
the effective conductance between $\GG'$ and $\GG''$
is finite; explicitly, we make $\sum_i 1/(\log k_i)<\infty$.
This constructs the graph $\GG$.
}\end{remark}
\subsection{Application to Random-Cluster Model}
\label{subsect: Application to Random-Cluster Model}

\begin{corollary}
{(Th.~\ref{thintro: beta-1 and pc,pu for Random-Cluster})}
\label{cor: `pc<pu' for G not in OHD and RC model}
Let $\GG$ be a unimodular transitive locally finite graph, not in $\OOO$.
	Fix the parameter $q\in [1,\infty)$.
	The gap between the left limit (when $p\nearrow p_c(q)$) and the right limit (when $p\searrow 
	p_u(q)$) of the expected degree of a base point $\origin$ with respect to the measure $\RC_{p,q}$ 
	satisfies:
	$$0<\beta_1(\GG)\leq \frac{1}{2} \Big( \RC_{p_u +, q} [\val (\origin)]-\RC_{p_c -, q} [\val 
	(\origin)]\Big).$$
\end{corollary}
Here, $\RC$ denotes either $\WRC$ or $\FRC$.

{\textsc{Proof}}: 
Consider the invariant coupling introduced by O.~H\"aggstr\"om, 
J.~Jonasson and R.~Lyons in \cite{HJL02a} of (all) the measures 
$\FRC_{p,q}$ and $\WRC_{p,q}$ (together) for $p\in[0,1]$ and $q\in [1,\infty)$
$$(X,\mu)\overset{\pi^{\RC}_{p,q}}{\rightarrow} (\{0,1\}^{\EE}, \RC_{p,q})$$
It provides two families of countable equivalence relations $\RRcluster_{p,q}$
(on the quotient space 
$\KK\backslash X$, once given a closed unimodular transitive group of automorphisms of $\GG$), 
one for $\FRC$ and one for $\WRC$.
The usefulness of that coupling is that it reflects the stochastic domination 
(see \cite[sect. 3]{HJL02a}); in particular, 
for a fixed parameter $q$ and $s<t$ (and denoting $\pi^{\RC}_{s,q}$ by $\pi_{s}$):
$$\mu\Bigl(\pi_{t}(e)=1 \mathrm{\ and\ } \pi_{s}(e) =0 \Bigr)\ \ =\ \ 
\mu(\pi_t(e)=1)-\mu(\pi_s(e)=1)$$
Take $s,t$ such that $s<p_c(q)\leq p_u(q)<t$, then Theorem~\ref{th: general pc<pu+ couplings} says that:
$$\displaystyle{\beta_{1}(\GG)\leq \frac{1}{2} \hskip-20pt
\sum_{{\begin{array}{c}
    \mbox{\footnotesize{edges\ }} e\\
    \mbox{\footnotesize{adjacent\ to\ }} \origin
    \end{array}}}\hskip-20pt 
\mu\bigl(\pi_t(e)=1)-\mu(\pi_s(e)=1\bigr)
}$$
Now, the right member is precisely:
$\frac{1}{2} \big( \RC_{t, q} [\val (\origin)]-\RC_{s, q} [\val (\origin)]\big).$ 
The monotonicity properties of the measures $\RC$ 
lead to the required inequality.
Indeed, monotonicity as well 
as left continuity of $p\mapsto \FRC_{p,q}[\val (\origin)]$
and right continuity of $p\mapsto \WRC_{p,q}[\val (\origin)]$
follow, like in \cite{HJL02b}, from the fact that $\FRC$ is an 
increasing (and $\WRC$ is a decreasing) limit of increasing (in $p$) continuous functions.
\endofproof


\section{Quasi-transitive graphs}
\label{sect: quasi-transitive graphs}

This section indicates how to extend the above results to
the context of quasi-transitive graphs, instead of just transitive ones.
There is no qualitative reversal, and just some quantitative adjustments.
The proofs are straightforward adaptions of those of the transitive case
with just slight changes of notation.
We first describe how to modify section~\ref{sect:Percolation on transitive graphs}.

Let $\GG=(\VV,\EE)$ be a locally finite quasi-transitive graph, 
$\HH$ a closed subgroup of $\mathrm{Aut}(\GG)$ whose action on
$\VV$ has finitely many orbits. Choose one vertex $\origin_1, \origin_2,
\cdots, \origin_{q}$ in each orbit and denote by $K_1=K_{\origin_1},K_2=K_{\origin_2}, \cdots,K_{q}=K_{\origin_{q}}$
its stabilizer.

\noindent
Let $(X,\mu)$ be a standard Borel probability space together with \\ 
	\indent $\bullet$ a probability measure-preserving (p.m.p.) action of 
    $\HH$, which is  \emph{essentially free}, and\\
    \indent $\bullet$  an $\HH$-equivariant Borel map $\pi:X\to \{0,1\}^{\EE}$.

\noindent
Divide out $X\times \GG$ by the diagonal action of 
$\HH$ to get the laminated space $\Lfull=\HH\backslash (X\times \GG)$: the 
\textbf{full lamination}.
Corresponding to the partition of $\VV$ into $\HH$-orbits $\VV_1,\VV_2, \cdots,\VV_{q}$, 
the transversal $X^{\bullet}:=\HH\backslash (X\times \VV)$
 identifies  with 
the disjoint union of standard Borel spaces
$K_1\backslash X \coprod K_2\backslash X \coprod\cdots \coprod K_{q}\backslash X$
via  $X\times \VV_i \ni  (hx,h\gamma \origin_i)
\mapsto K_i \gamma^{-1} x\in K_i\backslash X$.
The leaf of 
$\mu^{\bullet}$-almost every $x^{\bullet} \in  X^{\bullet}$ is isomorphic to 
$\GG$.

Define the \textbf{full equivalence relation} $\RRfull$ on $X^{\bullet}$ by 
$x^{\bullet}\RRfull y^{\bullet}$ if and only if $x^{\bullet}$ and 
$y^{\bullet}$ are vertices of the same $\Lfull$-leaf.

It inherits naturally an \textit{unoriented graphing}
and a \textit{smooth field of graphs}
(see  Section~\ref{sect: background Measured Equivalence 
    Relations}, and examples~\ref{examp: unoriented graphings}, 
\ref{examp: fields of graphs for transitive} of 
Section~\ref{sect:fields graphs, HD functions and L2 Betti numbers})
from the edge set $\HH\backslash (X\times \EE)$, where the graph 
associated with each point admits an 
isomorphism with $\GG$.

Consider, on $X^{\bullet}$, the  probability measure
$\mu^{\bullet}:=\frac{1}{T}( \frac{\Pi_{1*}\mu}{\Haar (K_1 )}
+\frac{\Pi_{2*}\mu}{\Haar (K_2 )}+\cdots +
\frac{\Pi_{p*}\mu}{\Haar (K_{q})}
) $, where $T=\frac{1}{\Haar (K_1)}
+\frac{1}{\Haar (K_2 )}+\cdots +
\frac{1}{\Haar (K_{q})}$ and $\Pi_{i*}\mu$ is the pushed-forward measure 
by $\Pi_{i}:X\to K_{i}\backslash X$. 
It is preserved by
the equivalence relation $\RRfull$
if and only if the group $\HH$ is unimodular (Th.~\ref{th: Rful preserves the measure iff unim.}).
Observe that with this choice (under the unimodularity asumption), the description depends neither on the choice of scaling of the Haar measure, nor on the choice of a particular orbit of vertices.
The {\bf first $L^2$ Betti number} of the equivalence relation $\RRfull$
is given by the formula:
$$\beta_{1}(\RRfull, \mu^{\bullet}):=\frac{1}{2T} 
    \sum_{i=1}^{q} \frac{1}{\Haar(K_i)} \sum_{j=1}^{n_{i}} 
    \scal{p(\mathbf{1}_{e_{i,j}})}{\mathbf{1}_{e_{i,j}}},
$$
where $p$ is the projection $p:C_{(2)}^{1}(\GG)\to\dzero 
\HD(\GG)\simeq \HD(\GG)/\Cmath$, and for each $i=1,2,\cdots,q$, the vectors
$\mathbf{1}_{e_{i,1}}, \mathbf{1}_{e_{i,2}}, \cdots, 
\mathbf{1}_{e_{i,n_{i}}}\in C_{(2)}^{1}(\GG)$ are the characteristic functions of (all) the
edges $e_{i,1}, e_{i,2}, \cdots, e_{i,n_i}$ adjacent to the orbit representative 
$\origin_i$ (see sect.~\ref{sect: some computation}).
It is not hard to check that this quantity doesn't depend on the particular unimodular quasi-transitive group of automorphisms $\HH$\footnote{this is
in fact true for the convex combination $\frac{1}{T} 
    \sum_{i=1}^{q} \frac{1}{\Haar(K_i)} \sum_{j=1}^{n_{i}} 
    \eta(e_{i,j})
$
associated with any $\mathrm{Aut}(\GG)$-invariant function  $\eta$ defined on the edges of $\GG$.}\,\footnote{It is most probably the case that the higher dimensional $L^2$ Betti numbers $\beta_{n}(\RRfull, \mu^{\bullet})$ are also
invariants of the graph (in fact, of the automorphism group of the unimodular graph, with a normalization given by an  appropriate combination of the Haar measures of the stabilizers of the vertices), but this would move us apart from the purpose of this paper.}.
One defines the {\bf first 
$\ell^2$ Betti number} of $\GG$ by the same formula (see def.~\ref{def: beta1 of the graph}):
$$\beta_1(\GG):=\frac{1}{2T} 
    \sum_{i=1}^{q} \frac{1}{\Haar(K_i)} \sum_{j=1}^{n_{i}} 
    \scal{p(\mathbf{1}_{e_{i,j}})}{\mathbf{1}_{e_{i,j}}}.
$$
Clearly, $\beta_1(\GG)=0$ if and only if $\GG$ belongs to $\OOO$
(see footnote~\ref{footn: beta1 and OHD}).
\begin{example}
The leading example is the group $\HH=\Zmath/r\Zmath\ast \Zmath/s\Zmath$ 
acting on its Bass-Serre tree
$\GG$: the bipartite tree with valencies $r$ and $s$.
The stabilizers are $K_1=\Zmath/r\Zmath$ and $K_2=\Zmath/s\Zmath$
and the components of $X^{\bullet}\simeq (\Zmath/r\Zmath)\backslash X\coprod
(\Zmath/s\Zmath)\backslash X$ identifie with pieces of $X$ of measure
$\frac{1}{r}$ and $\frac{1}{s}$. By considering the cost of the graphing
inherited by
$\RRfull$ (see example~\ref{examp: fields of graphs for transitive} of 
Section~\ref{sect:fields graphs, HD functions and L2 Betti numbers} and Section~\ref{sect: background Measured Equivalence Relations}),
one computes:
$\beta_1(\GG)=\beta_1(\RRfull, \mu^{\bullet})=(\frac{1}{r}+\frac{1}{s})^{-1}[1-(\frac{1}{r}+\frac{1}{s})]$. 
In the particular case of the regular tree of even valency $r=s=2t$, 
$\beta_1(\RRfull, \mu^{\bullet})=t-1=\beta_1(\FF_{t})$ coincide for $\HH=\Zmath/2t\Zmath\ast \Zmath/2t\Zmath$ or $\HH=\FF_{t}$.
\end{example}

Thanks to the map $\pi:X\to \{0,1\}^{\EE}$,  each leaf of $\Lfull$ becomes a colored graph.
The \textbf{cluster lamination} is obtained by removing
all the $0$-colored edges.
Define the \textbf{cluster equivalence relation} $\RRcluster$ on $X^{\bullet}$ by 
$x^{\bullet}\RRcluster y^{\bullet}$ if and only if $x^{\bullet}$ and 
$y^{\bullet}$ are vertices of the same $\Lcluster$-leaf. It is a 
subrelation of $\RRfull$.

The proof of Theorem~\ref{th: virt select clusters in/out OHD} extends
to quasi-transitive graphs with no modification.

\medskip
\noindent
\textbf{Theorem~\ref{th: virt select clusters in/out OHD}*} \hskip5pt
{\em Assume that the closed subgroup $\HH$ is unimodular and $\pi_*\mu$ admits an $\HH$-equivariant 
virtually selected cluster. 
    If $\GG$ belongs to $\OOO$, then $\mu$-a.e.\ virtually selected
    cluster belongs to $\OOO$.
    If $\GG$ doesn't belong to $\OOO$, then $\mu$-a.e.\ virtually
    selected cluster doesn't belong to $\OOO$.\\
    In particular, on the Borel set of subgraphs with one infinite cluster,
    the infinite cluster belongs (resp. doesn't belong) to $\OOO$ a.s. i{f}{f} $\GG$ belongs (resp. doesn't belong) to $\OOO$.}

\bigskip

We now turn to give the modifications in the statements of 
the quantitative estimates of section~\ref{sect: Non uniq. pahse Harm. Dir. Funct}.

\medskip
\noindent
\textbf{Theorem~\ref{th: general pc<pu+ couplings}*} \hskip5pt
{\em 
Consider a unimodular quasi-transitive group $\HH$ of automorphisms of $\GG$, 
two $\HH$-invariant percolations $\mu_{1}$ and $\mu_{2}$ on $\GG$
and an 
\textbf{$\HH$-equivariant coupling}
$$
\begin{array}{rcl}
    & (X,\mu) &\\
    {}^{\pi_{1}}\swarrow & & \searrow^{\pi_{2}}\\
    (\{0,1\}^{\EE},\mu_{1}) & & (\{0,1\}^{\EE}, \mu_{2})
\end{array}
$$
\vskip-10pt
Assume that \\ 
	\indent 1. $\mu_{1}$-a.e. cluster belongs to $\OOO$,\\
	\indent 2. $\pi_2$ has an $\HH$-equivariant selected cluster defined on a 
 	non-null set\,\footnote{For instance, if
 	$\mu_{2}$ has a non-null set of subgraphs with exactly one 
   infinite cluster.},\\
   then
$$\displaystyle{\beta_{1}(\GG)\leq {\frac{1}{2T}}\Bigl[
\sum_{i=1}^{q}
\frac{1}{\Haar (K_i) }
\hskip-20pt
\sum_{{\begin{array}{c}
    \mbox{\footnotesize{\em edges\ }} e\\
    \mbox{\footnotesize{\em adjacent\ to\ }} \origin_i
    \end{array}}}\hskip-20pt
\mu\Bigl(\pi_{2}(e)=1 \mathrm{\ and\ } \pi_{1}(e) =0 \Bigr)
}
\Bigr]
.$$
}
The only modifications in the proof of Theorem~\ref{th: general pc<pu+ couplings} are the bound on $\mathrm{cost}(\Psi_{2})$:
$$\mathrm{cost}(\Psi_{2})\leq \frac{1}{2 \mu^{\bullet}(Y)}
\frac{1}{T}
\sum_{i=1}^{q} \frac{1}{\Haar(K_i)}\hskip-20pt
\sum_{{\begin{array}{c}
    \mbox{\footnotesize{edges\ }} e_{i,j}\\
    \mbox{\footnotesize{adjacent\ to\ }} \origin_i
    \end{array}}}
\hskip-20pt
\mu\Bigl(\pi_{2}(e)=1 \mathrm{\ and\ } \pi_{1}(e) =0 \Bigr),$$
 and the definition of the measurable subset $Z\subset X^{\bullet}$:
For at least one of the $\origin_i$, the set of those $x\in X$ whose selected cluster contains $\origin_i$ is non-null. Take for $Z$ its image in $X^{\bullet}$.
\endofproof

\bigskip
\noindent
\textbf{Application to Bernoulli Percolation} (subsection~\ref{subsect: Application to Bernoulli Percolation}*).

\medskip
\noindent
\textbf{Corollary~\ref{cor: pc<pu for G not in OHD}*} 
\hskip5pt
{\em 
	Let $\GG$ be a unimodular quasi-transitive locally finite graph.
    If $\GG$ doesn't belong to $\OOO$, then the nonuniqueness phase 
    interval of Bernoulli percolation has non-empty interior:
    $$p_{c}(\GG) < p_{u}(\GG)$$
	More precisely, 
	$$0<\beta_{1}(\GG)\leq \frac{1}{2T} \sum_{i=1}^{q} \frac{\val(\origin_i)}{\Haar (K_i) } \bigl(p_{u}(\GG)-p_{c}(\GG)\bigr).$$
where $\val(\origin_i)$ is the number of edges in $\GG$ that are adjacent to $\origin_i$.}\endofproof

\bigskip
\noindent
\textbf{Application to Random-Cluster Model} (subsection~\ref{subsect: Application to Random-Cluster Model}*).

\medskip
\noindent
\textbf{Corollary~\ref{cor: `pc<pu' for G not in OHD and RC model}*} 
\hskip5pt
{\em 
Let $\GG$ be a unimodular quasi-transitive locally finite graph, not in $\OOO$.
	Fix the parameter $q\in [1,\infty)$.
	The gap between the left limit (when $p\nearrow p_c(q)$) and the right limit (when $p\searrow p_u(q)$) of the expected degree of a base point $\origin$ with respect to the measure $\RC_{p,q}$ 
	satisfies:
	$$0<\beta_1(\GG)\leq \frac{1}{2T} \sum_{i=1}^{q} \frac{1}{\Haar (K_i )}
	\Big( \RC_{p_u +, q} [\val (\origin_i)]-\RC_{p_c -, q} [\val 
	(\origin_i)]\Big).$$
}\\
Here, $\RC$ either denotes $\WRC$ or $\FRC$, and 
$\RC_{p,q}[\val(\origin_i)]$ is the mean degree 
of the vertex $\origin_i$ in the random subgraph
for the random-cluster measure $\RC$ with parameters 
$p$ and $q$. Also $ \RC_{p_u +, q} [\val (\origin_i)]:=\lim_{p\searrow p_u(q)}\RC_{p,q}[\val(\origin_i)]$ and 
$\RC_{p_c -, q} [\val 
	(\origin_i)]:=\lim_{p\nearrow p_c(q)}\RC_{p,q}[\val(\origin_i)]$.
\endofproof

\section{Harmonic Dirichlet Functions and $\ell^2$ Cohomology}
\label{sect:Harmonic Dirichlet Functions and L2 Cohomology}

The main result of the section
is the technical Proposition~\ref{prop:cap H(S-s) to H(S-t) isomorphic with 
HD for simpl-compl} that will allow to make the connection 
between harmonic Dirichlet functions and  
the definitions of $L^2$ Betti numbers for equivalence relations in \cite{Gab02}
(see the proof of Theorem~\ref{thm: Betti 1 and graphings}).
However, as a leading and motivating example, we will consider 
the following well-known result 
relating harmonic Dirichlet functions
with the first $\ell^2$ Betti number.
\begin{theorem}\label{thm: beta-1 and HD of a Cayley graph}
    Let $\GGamma$ be a finitely generated group.
    Its first $\ell^{2}$ Betti number $\beta_{1}(\GGamma)$ is not zero
    if and only if any of its Cayley 
    graphs admits nonconstant harmonic Dirichlet functions.
\end{theorem}

\subsection{Harmonic Dirichlet Functions\ldots}
Let $\GG=(\VV,\EE)$ be a connected graph\,\footnote{In the whole 
Section~\ref{sect:Harmonic Dirichlet Functions and L2 Cohomology}, except in the examples, $\GG$ is not assumed 
to have any kind of symmetries or automorphisms.} with bounded degree. 
The tail and head of an oriented edge $\hat{e}$ are denoted by $\hat{e}^{-}$ and $\hat{e}^{+}$.

Denote by ${\cal F}(\VV)$ or $C^{0}(\GG)$ the space of all complex-valued 
functions ($0$-cochains) on $\VV$
and by $C^{1}(\GG)$ the space of $1$-cochains\,\footnote{{i.e.\ } the 
space of anti-symmetric functions on the set of oriented edges:
$f(\check{e})=-f(\hat{e})$ where $\check{e}$ is the edge $\hat{e}$ with the reverse orientation}.
Define the coboundary, ``boundary'' and Laplace 
maps:\\
$\begin{array}{cccccrl}
    \dzero : & C^{0}(\GG) & \longrightarrow & C^{1}(\GG) & & 
    \dzero f(\hat{e}) =& \displaystyle{f(\hat{e}^{+})-f(\hat{e}^{-})} \phantom{ \displaystyle{\sum_{\hat{e}: \hat{e}^{+}=v} g(\hat{e})}}\\
    \dzero^{*}: & C^{1}(\GG) & \longrightarrow & C^{0}(\GG) & & 
    \dzero^{*} g(v) =&
    \displaystyle{\sum_{\{\hat{e} \colon \hat{e}^{+}=v\}} g(\hat{e})}\\
    \Delta: & C^{0}(\GG) & \longrightarrow & C^{0}(\GG) & & 
    \Delta f (v)=\dzero^{*} \dzero f (v)  =&
    \displaystyle{\sum_{\{\hat{e} \colon \hat{e}^{+}=v\}}\bigl(f(\hat{e}^{+})-f(\hat{e}^{-})\bigr)}\\
    &&&&&=& \displaystyle{\val(v) f(v)-\sum_{\{\hat{e} \colon \hat{e}^{+}=v\}}f(\hat{e}^{-})}
\end{array}$\\
The spaces of $\ell^{2}$-cochains will be denoted 
by $C^{0}_{(2)}(\GG)$ and $C^{1}_{(2)}(\GG)$. 
By definition, the space of \textbf{harmonic Dirichlet functions} on $\GG$ 
is the space of functions whose value at each $v$ equals the mean of 
the values at its neighbors ($\Delta(f)=0$, i.e.\  $f$ is harmonic) and with 
 coboundary in $\ell^2$ ($f$ has finite energy or finite Dirichlet sum):
$$\HD(\GG):=\{f\in {\cal F}(\VV) : \dzero f\in 
C^{1}_{(2)}(\GG)\textrm{ and } \Delta f = 0 \}.$$
The kernel of $\dzero$ clearly consists of
the constant functions (since 
$\GG$ is connected), so that naturally 
\begin{eqnarray}
    \HD(\GG)/ \Cmath\ \simeq\ \dzero(\HD(\GG))=\mathrm{Im}\dzero\cap 
    C_{(2)}^{1}(\GG)\cap {\mathrm{Ker}}\dzero^{*}.\label{eqn: HD(G) simeq Im cap C1 cap ker}
\end{eqnarray}

\medskip

As for the $\ell^2$ cohomology, it is not really that of the graph 
$\GG$ that is of interest, since, for example, for Cayley graphs it is 
too sensitive to changes of generators (think of Cayley graphs of 
$\Zmath$, where $\HHH=0$ or $\not=0$ according to whether the 
generating system is $(1)$ or $(2,3)$). One has first to ``fill in the 
holes'' of the graph: \\
Consider a  simply-connected $2$-dimensional complex $\Sigma$, 
with $\GG$ as $1$-skeleton.

\begin{example}
    For instance, let $\GGamma$ be a group given by a presentation with $g$ generators 
and $r$ relators. Recall that the \textbf{Cayley complex} of the 
presentation is a $2$-dimensional complex built from a bouquet of $g$ oriented circles labeled 
by the generators, together with $r$ oriented disks labeled by the 
relators glued along their boundary to this $1$-dimensional skeleton 
by following successively the circles associated with the labeling 
relator. Its fundamental group is (isomorphic to) $\GGamma$.
The universal cover $\Sigma$ of the Cayley complex is a (simply 
connected) $2$-dimensional complex (with a free action of $\GGamma$ and) with the 
Cayley graph of $\GGamma$ as $1$-skeleton. 
\end{example}

More generally, $\Sigma$ can be obtained from $\GG$ by gluing one oriented disk 
(thought of as a polygon) along its boundary to each 
circuit (and the opposite orientation for the reverse circuit).

\medskip

Denote by $C^{2}(\Sigma)$ the space of $2$-cochains, {i.e.\ } anti-symmetric functions 
on the oriented $2$-cells --disks--), by $C_{(2)}^{2}(\Sigma)$ 
the space of those that are $\ell^2$ (i.e.\  $\sum_{\sigma} 
h(\sigma)^2<\infty$, where the sum is over all the $2$-cells $\sigma$).
The boundary of a $2$-cell $\sigma$ in $\Sigma$ being a $1$-cycle $\dun^{*} 
\sigma$, one defines the \emph{coboundary} $\dun$ by $\dun g(\sigma)=g (\dun^{*} 
\sigma)$. 
$$\begin{array}{ccccc}
C^{0}(\Sigma) &\overset{\dzero}{\longrightarrow} &
C^{1}(\Sigma) &\overset{\dun}{\longrightarrow} &
C^{2}(\Sigma)
\end{array}$$
By taking the adjoint, $\smash{C_{(2)}^{1}(\Sigma)\cap 
{\mathrm{Ker}}\dzero^{*}=\bigl[{\dzero}C_{(2)}^{0}(\Sigma)\bigr]^{\perp 
C_{(2)}^{1}}}$, the orthogonal of ${\mathrm{Im}}\dzero_{(2)}={\dzero}C_{(2)}^{0}(\Sigma)$ in $C_{(2)}^{1}(\Sigma)$.
Since $\Sigma$ is simply connected,
$\mathrm{Im}\dzero=\mathrm{Ker}\dun$. 
From formula (\ref{eqn: HD(G) simeq Im cap C1 cap ker}) we get the natural 
isomorphisms
\begin{eqnarray}
    \HD(\GG)/ \Cmath\ \simeq\ \mathrm{Ker}\dun\cap 
    C_{(2)}^{1}(\Sigma) \cap 
    \bigl[{\dzero}C_{(2)}^{0}(\Sigma)\bigr]^{\perp}\simeq \frac{\mathrm{Ker}\dun\cap
    C_{(2)}^{1}(\Sigma)} {{\overline{\mathrm{Im}}\dzero_{(2)}}},
    \label{eqn: HD isom Harmonic fct in C-1-(2)}
\end{eqnarray}
where ${\overline{\mathrm{Im}}\dzero_{(2)}}$ is the closure of the
space ${\mathrm{Im}}\dzero_{(2)}$.

\subsection{\ldots and $\ell^2$ cohomology}
To define $\ell^2$ cohomology of $\Sigma$, one is led to consider
$\ell^2$ cochains and restrictions of the coboundary maps:
$$
C_{(2)}^{0}(\Sigma) \overset{\dzero_{(2)}}{\longrightarrow}
C_{(2)}^{1}(\Sigma) \overset{\dun_{(2)}}{\longrightarrow}
C_{(2)}(\Sigma)$$

Say that $\Sigma$ is \textbf{uniformly locally bounded} (ULB) if it
admits a uniform bound $M$ s.t. each vertex (resp.\   edge) belongs to
at most $M$ edges (resp.\   $2$-cells), and the boundary of each
$2$-cell has length at most $M$.  In this situation, $\dun_{(2)}$ is a
bounded operator, and the standard \textbf{first reduced $\ell^2$
cohomology space} of $\Sigma$ is defined as the Hilbert space $$\HHH
(\Sigma):=
\smash{\frac{\mathrm{Ker}\dun_{(2)}}{\overline{\mathrm{Im}}\dzero_{(2)}}}=\frac{\mathrm{Ker}\dun\cap
    C_{(2)}^{1}(\Sigma)} {{\overline{\mathrm{Im}}\dzero_{(2)}}}.$$ 
It follows from (\ref{eqn: HD isom Harmonic fct in C-1-(2)}) that for 
a ULB $\Sigma$
\begin
{eqnarray}
    \HD(\GG)/ \Cmath\ \simeq \HHH
(\Sigma).\label{eq: HD isom l2 cohomology}
\end{eqnarray}

\begin{example} $\Sigma$ is ULB if it comes from the Cayley complex
    of a finitely presented group.  One can show that any other
    simply connected, \emph{cocompact} free $\GGamma$-complex $\Sigma'$
    leads to a ($\GGamma$-equivariantly) isomorphic $\HHH(\Sigma')$, so
    that the non-triviality of $\HHH(\Sigma)$ is an invariant of the
    group $\GGamma$, and not only of the complex $\Sigma$.  The first
    $\ell^2$ Betti number $\beta_{1}(\GGamma)$ is the von Neumann
    $\GGamma$-dimension $\dim_{\GGamma} \HHH (\Sigma)$ and we retain from
    this dimension theory that it vanishes i{f}{f} $\HHH (\Sigma)=\{0\}$.
    
    Observe that some finiteness condition on $\Sigma$ is however
    necessary since the free group $\FF_{2}$, for example, acts on its
    Cayley tree as well as on the product of the tree with a line, whose
    $1$-skeleton admits (resp.\   doesn't admit) nonconstant harmonic Dirichlet
    functions.  In case the presentation of $\GGamma$ is not finite
    ($r=\infty$), $\Sigma$ is no longer $\GGamma$-cocompact, and
    $\dun_{(2)}$ is no longer a continuous (=bounded) map.
\end{example}

The general way to proceed to define $\ell^2$ cohomology for a complex that is
not ULB, in the spirit of J.~Cheeger and M.~Gromov \cite{CG86},
consists in approximating $\Sigma$ by its ULB subcomplexes.  Consider
the directed set of ULB subcomplexes $\Sigma_{t}$ of $\Sigma$,
directed by inclusion and the inverse system of reduced $\ell^2$
cohomology spaces $H_{(2)}^{n}(\Sigma_{t})$ of $\Sigma_{t}$ with the
maps $H_{(2)}^{n}(\Sigma_{s})\to H_{(2)}^{n}(\Sigma_{t})$ induced by
inclusion $\Sigma_{s}\supset \Sigma_{t}$ (denoted by $s\geq t$). 
Then define the reduced $\ell^2$ cohomology as the inverse limit
$\displaystyle{H_{(2)}^{n}(\Sigma):=\lim_{\leftarrow}}
H_{(2)}^{n}(\Sigma_{t})$.

\medskip

In our context, all the ULB complexes $\Sigma_{t}$, as well as 
$\Sigma$ itself, share the same ULB $1$-skeleton $\GG$.
Thus, the first reduced $\ell^2$ cohomology spaces $\HHH(\Sigma_{t})$
are each the quotient of the subspace $\mathrm{Ker} \Bigl(C_{(2)}^{1}(\GG)\overset{\dun}{\longrightarrow} 
C_{(2)}^{2}(\Sigma_{t})\Bigr)$ of $C_{(2)}^{1}(\GG)$ (organized into 
an inverse system by inclusion), by the common subspace ${\overline{\mathrm{Im}} \dzero_{(2)}}$.
It follows that:
\begin{eqnarray}
    \displaystyle{\HHH(\Sigma)=\lim_{\leftarrow }} \HHH(\Sigma_{t})
    =\frac{\bigcap_{t}\mathrm{Ker} 
\bigl(C_{(2)}^{1}(\GG)\overset{\dun}{\longrightarrow} 
C_{(2)}^{2}(\Sigma_{t})\bigr)}{\overline{\mathrm{Im}} 
\dzero_{(2)}}=
\frac{\mathrm{Ker}\dun\cap C_{(2)}^{1}(\GG)}{\overline{\mathrm{Im}} 
\dzero_{(2)}} \label{eqn: H(2) as limit of inverse system}
\end{eqnarray}
and by formula (\ref{eqn: HD isom Harmonic fct in C-1-(2)}), valid for any 
$\Sigma$:
\begin{eqnarray}
    \HHH(\Sigma)\ \simeq\ \HD(\GG)/ \Cmath.
\end{eqnarray}

\begin{example} Let $\GGamma$ be finitely generated, but not 
    necessarily finitely presented.
    For a $\GGamma$-complex $\Sigma$, the $\Sigma_{t}$ are moreover required to
	be $\GGamma$-invariant (and cocompact) and the $\ell^2$ Betti numbers of the 
	$\GGamma$-action on $\Sigma$ are 
    defined by keeping track of the $\GGamma$-dimensions:
    $$
    \begin{array}{rcl}
	\beta_{n}(\Sigma, \GGamma)&:=& \sup_{t} \dim_{\GGamma} \overline{\mathrm{Im}}(\bar{H}_{(2)}^{n}(\Sigma)\to 
	\bar{H}_{(2)}^{n}(\Sigma_{t})) \\
	&=&
	\sup_{t}
	\dim_{\GGamma} \cap_{s\geq t} \overline{\mathrm{Im}}(\bar{H}_{(2)}^{n}(\Sigma_{s})\to 
	\bar{H}_{(2)}^{n}(\Sigma_{t}))
    \end{array}$$
    For simply connected $\GGamma$-complexes $\Sigma$, 
    the value $\beta_{1}(\Sigma, \GGamma)$ doesn't depend on a 
    particular choice of $\Sigma$, so that for $\Sigma$ constructed from a Cayley complex of $\GGamma$, 
    $$\beta_{1}(\GGamma)=\beta_{1}(\Sigma,\GGamma)= \dim_{\GGamma} \HD(\GG)/ 
    \Cmath$$ vanishes if and only if $\HD(\GG)/ 
    \Cmath=\{0\}$. This proves Theorem~\ref{thm: beta-1 and HD of a Cayley graph}.
\end{example}

\bigskip
\noindent
Let's quote for further use the observation that the space of 
formula~(\ref{eqn: H(2) as limit of inverse system}) may be obtained by considering a exhausting sequence 
instead of the whole inverse system:

\begin{proposition}
	Let $\GG$ be a graph with finite degree, $\Sigma$ a 
	simply-connected $2$-dimensional complex with {$1$-skeleton~$\GG$}.
    If $(\Sigma_{t})_{t\in \Nmath}$ is an increasing and exhausting 
	sequence of ULB 
    subcomplexes of $\Sigma$, then for any fixed $t$
    $$\cap_{s\geq t} {\mathrm{Im}}(\HHH(\Sigma_{s})\to 
    \HHH(\Sigma_{t}))=\frac{\mathrm{Ker}\dun\cap 
	C_{(2)}^{1}(\GG)}{\overline{\mathrm{Im}}
    \dzero_{(2)}}$$ doesn't depend on $t$ and is NATURALLY isomorphic with 
    $\HD(\GG)/ 
    \Cmath$.
\end{proposition}
The connection with the simplicial framework of \cite{Gab02} 
is made by considering a double barycentric subdivision $\Sigma^{*}$ 
of $\Sigma$, with the exhaustion $\Sigma^{*}_{t}$ corresponding to 
the subdivision of $\Sigma_{t}$. Since for each $t$, 
$\HHH(\Sigma_{t})$ and $ \HHH(\Sigma^{*})$ are naturally isomorphic, it follows that 
\begin{proposition} \label{prop:cap H(S-s) to H(S-t) isomorphic with 
    HD for simpl-compl}
    For any fixed $t$, $\cap_{s\geq t} {\mathrm{Im}}(\HHH(\Sigma_{s}^{*})\to 
    \HHH(\Sigma_{t}^{*}))$ doesn't depend on $t$ and is NATURALLY isomorphic with $\dzero(\HD(\GG))$ and $\HD(\GG)/ 
    \Cmath$.
\end{proposition}

\section{Fields of Graphs, Harmonic Dirichlet Functions and $L^{2}$~Betti 
Numbers for Equivalence Relations}
\label{sect:fields graphs, HD functions and L2 Betti numbers}

Let $(X,\mu)$ be a standard Borel space with a probability 
measure $\mu$ and $\RR$ a measure-preserving Borel equivalence relation
with countable classes.

Recall from \cite{Gab02} that an \textbf{$\RR$-equivariant field $x\mapsto \Sigma_{x}$ 
of simplicial complexes} is a measurable 
assignment to each $x\in X$ of a simplicial complex $\Sigma_{x}$,
together with an ``action'' of $\RR$, i.e.\  
with the measurable data of a simplicial isomorphism, for every 
$(x,y)\in \RR$,
$\psi_{x,y}:\Sigma_{y}\to \Sigma_{x}$ such that 
$\psi_{x,y}\psi_{y,z}=\psi_{x,z}$ and $\psi_{z,z}=id_{\Sigma_{z}}$.
It is \textbf{smooth} if the action on the vertices admits a Borel 
fundamental domain.
It is \textbf{smooth uniformly locally bounded } 
if there is a uniform bound $N$ on the degree of the $1$-skeleton of the
$\Sigma_{x}$, and there is a Borel fundamental domain that meets 
each $\Sigma_{x}$ in at most $N$ vertices.

\begin{example}\label{examp: unoriented graphings}
    Let $\RR$ be a p.m.p.\ countable Borel equivalence relation on the probability standard Borel
    space $(X,\mu)$. 
    An {unoriented graphing} $\Psi$ over $\RR$ (see Section~\ref{sect: background Measured Equivalence 
    Relations}) defines an
    {$\RR$-equivariant field of graphs} $x\mapsto 
    \Psi_{x}$ with vertex set $\RR$ itself, which is smooth.\\
    - The vertex set of $\Psi_{x}$ is the set $\{(x,y)\in \RR\}$, i.e.\  the 
    set of elements of $\RR$ with first coordinate $x$.\\
    - Two vertices $(x,y)$ and $(x,z)$ of $\Psi_{x}$ are neighbors if and only 
    if $(y,z)$ belongs to $\Psi$, i.e.\  i{f}{f} the second coordinates are 
    neighbor for $\Psi$\\
    - The left action of $\RR$ on itself $(w,x).(x,y)=(w,y)$ and thus on the 
    set of vertices induces a natural action on the field:
    $(w,x) :\left(
    \begin{array}{ccc}
	\Psi_{x} & \longrightarrow & \Psi_{w} \\
	\left[(x,y),(x,z)\right] & \mapsto & \left[(w,y),(w,z)\right]
    \end{array}\right)
    $.\\
    - The ``diagonal'' set  $\{(x,x):x\in X\}$ of vertices forms a Borel 
    fundamental domain.
\end{example}
This example contains as main applications the various equivariant fields 
of graphs (described below) relevant for percolation theory.
\begin{example}\label{examp: fields of graphs for Cayley}
    In the context of Section~\ref{sect:Percolation on Cayley Graphs} for 
    a Cayley graph $\GG$, 
    \begin{description}
	\item $(\RRfull,\  x\mapsto \GG)$\hskip 5pt
	The full 
	lamination $\Lfull$ (Section~ \ref{subsect: full equiv. 
	rel. for Cayley graphs}) defines an unoriented graphing over 
	$\RRfull$ (see Section~\ref{sect: background Measured Equivalence 
	Relations}). In the corresponding smooth $\RRfull$-equivariant field 
	$x\mapsto \Lfull_{x}$, 
	each $\Lfull_{x}$ admits a canonical isomorphism with $\GG$.
	\item $(\RRcluster, \ x\mapsto \pi(x)(\origin))$\hskip 5pt
	The cluster lamination $\Lcluster$ defines an unoriented 
	graphing over $\RRcluster$. In the corresponding smooth 
	$\RRcluster$-equivariant field 
	$x\mapsto \Lcluster_{x}$, each $\Lfull_{x}$ is isomorphic to the cluster 
	$\pi(x)(\origin)$ of $\origin$ in $\pi(x)$.
	\item $(\RRfull,\ x\mapsto \pi(x))$ \hskip 5pt
	The cluster lamination $\Lcluster$ defines also an unoriented 
	graphing over $\RRfull$ and in the corresponding field 
	$x\mapsto \Psi_{x}$, 
	each $\Psi_{x}$ is isomorphic to the subgraph 
	$\pi(x)$, which is non-connected in general.
    \end{description}
\end{example}

\begin{example}\label{examp: fields of graphs for transitive}
    In the context of Section~\ref{sect:Percolation on transitive 
    graphs} for a transitive, locally finite graph $\GG$, 
    \begin{description}
	\item $(\RRfull,\ x^{\bullet}\mapsto \GG)$ \hskip 5pt The full
	lamination $\Lfull$ (Section~ \ref{subsect: full equi rel for
	transitive graph}) defines an unoriented graphing over
	$\RRfull$.  In the corresponding smooth $\RRfull$-equivariant
	field $x^{\bullet}\mapsto \Lfull_{x^{\bullet}}$ each
	$\Lfull_{x^{\bullet}}$ admits a (non-canonical) isomorphism
	with $\GG$.  However, each representative $x\in X$ of
	$x^{\bullet}$ defines an isomorphism
	$j_{x}:\Lfull_{x^{\bullet}}\simeq \GG$ and for another one
	$y=\kk x$, $j_{y}=\kk j_{x}$, where $\kk\in \KK$ so that this
	isomorphism is canonical up to an element of $\KK$.  \item
	$(\RRcluster, \ x^{\bullet}\mapsto \pi(x)(\origin))$ \hskip 5pt The
	cluster lamination $\Lcluster$ defines an unoriented graphing
	over $\RRcluster$.  The corresponding field
	$x^{\bullet}\mapsto \Lcluster_{x^{\bullet}}$ assigning to
	$x^{\bullet}$ its leaf (graph) in the lamination $\Lcluster$
	is a smooth uniformly locally bounded $\RRcluster$-equivariant
	field of connected graphs.  Each $\Lfull_{x^{\bullet}}$ is
	isomorphic to the cluster $\pi(x)(\origin)$ of $\origin$ for a (any)
	representative $x\in X$ of $x^{\bullet}$.  Two representatives
	give subgraphs of $\GG$ that are isomorphic under an element
	of $\KK$.  The graph $\Lcluster_{x^{\bullet}}$ belongs to
	$\OOO$ i{f}{f} $\pi(x)(\origin)$ belongs to $\OOO$ for any
	representative $x$ of $x^{\bullet}$.
    \end{description}
\end{example} 

\begin{example}\label{examp: smooth fields of graphs+restrictions}
    (Restrictions) If $x\mapsto \Psi_{x}$ is a {smooth
    $\RR$-equivariant field of graphs} and $Y$ is a Borel subset then
    restricted to $Y$, the fields $Y\ni x \mapsto \Psi_{x}$ is a {smooth 
	$\RR_{\vert Y}$-equivariant field of graphs}, where $\RR_{\vert Y}$ is
	the restriction of $\RR$ to $Y$.
\end{example}

\medskip Also recall from \cite{Gab02} that there is a well-defined
notion of $L^{2}$ Betti numbers $\beta_{n}(\RR,\mu)$ for a measure-preserving Borel
equivalence relation $\RR$ with countable classes, which uses the
notion of equivariant fields of simplicial complexes and the
\textit{von Neumann dimension} $\dim_{\RR}$ associated with the von
Neumann algebra of the equivalence relation and 
the measure~$\mu$.

\begin{theorem}\label{thm: Betti 1 and graphings}
    Let $\RR$ be a measure-preserving equivalence relation with countable 
    classes on the standard 
    Borel probability measure space $(X,\mu)$. Consider 
    a smooth uniformly locally bounded $\RR$-equiva\-riant field $x\mapsto 
    \GG_{x}$ of connected graphs.
    Then 
    $$\beta_{1}(\RR,\mu)=\dim_{\RR}\int^{\oplus}_{X} \dzero(\HD(\GG_{x})) \ d\mu(x)=\dim_{\RR}\int^{\oplus}_{X} \HD(\GG_{x})/\Cmath \ d\mu(x).$$
\end{theorem}
\noindent
Since $\dim_{\RR} H=0$ if and only if $H=\{0\}$, one gets:
\begin{corollary}\label{cor: relate beta-1 and HD functions for 
    relations}
    For a smooth uniformly locally bounded $\RR$-equivariant field $x\mapsto 
    \GG_{x}$ of connected graphs, 
    $\beta_{1}(\RR,\mu)=0$ if and only if $\mu$ a.s.\ 
    ${\GG_{x}}\in\OOO$.
\end{corollary}
\begin{remark}\label{rem: unoriented vs oriented graphings}
    Since a generating (oriented) graphing in the sense of 
    \cite{Lev95,Gab00} (see also Section~\ref{sect: background Measured Equivalence 
    Relations}) defines 
    an unoriented graphing and thus a smooth $\RR$-equivariant field of 
    connected graphs, the above Corollary~\ref{cor: relate beta-1 and HD functions for 
    relations} reduces to Theorem~\ref{thmintro: Betti 1 and graphings} of the 
    introduction.
\end{remark}

\proof 
By Theorem/Definition \cite[Th.~3.13, D\'ef.~3.14]{Gab02}, $\beta_{1}(\RR,\mu)$ is the 
first $L^{2}$ Betti number $\beta_{1}(\Sigma,\RR,\mu)$ of 
ANY smooth  $\RR$-equivariant field of 
\emph{simply connected} ($2$-dimensional, say) simplicial 
complexes $\Sigma$. 
It can be computed \cite[prop.~3.9]{Gab02} by using any 
exhausting increasing sequence $(\Sigma_{s})_{s\in \Nmath}$ of
$\RR$-invariant uniformly locally bounded (ULB) sub-complexes 
by the following formula~(\ref{eq:prop. 3.9 of Gab02}):
\begin{eqnarray}
   \beta_{1}(\Sigma,\RR,\mu)  & = & \lim_{s\to \infty} \!\!\nearrow
   \lim_{s\geq t\ s\to \infty}\!\!\!\!\!\searrow
   \dim_{\RR}\  \overline{\mathrm{Im}}\bigl[\bar{H}_{1}^{(2)}(\Sigma_{t}, \RR,\mu)\to 
   \bar{H}_{1}^{(2)}(\Sigma_{s}, \RR,\mu)\bigr]\label{eq:prop. 3.9 of Gab02}\\
     & = & \lim_{t\to \infty} \!\!\nearrow
   \lim_{s\geq t\ s\to \infty}\!\!\!\!\!\searrow
   \dim_{\RR}\  \overline{\mathrm{Im}}\bigl[\HHH(\Sigma_{s}, \RR,\mu)\to 
   \HHH(\Sigma_{t}, \RR,\mu)\bigr]\label{eq:rank of image equal for f 
   and its adjoint}\\
    & = & \lim_{t\to \infty} \!\!\nearrow \dim_{\RR}\!\!\!\bigcap_{s\geq 
    t\ s\to \infty}\!\!\!
    \overline{\mathrm{Im}}\bigl[\HHH(\Sigma_{s}, \RR,\mu)\to \HHH(\Sigma_{t}, 
    \RR,\mu)\bigr]\label{eq:exchanging dim and lim}
\end{eqnarray}

The equality~(\ref{eq:rank of image equal for f and its adjoint}) holds 
by duality between homology and cohomology
because, just as in usual linear algebra, taking dual does 
not alter the dimension of the image. Equality (\ref{eq:exchanging dim and 
lim}) is due to the continuity of dimension, since $\overline{\mathrm{Im}}\bigl[\HHH(\Sigma_{s}, \RR,\mu)\to \HHH(\Sigma_{t}, 
\RR,\mu)\bigr]\subset \HHH(\Sigma_{t}, 
\RR,\mu)$ decreases with~$s$.
Now, for a fixed $t$, one has the Hilbert integral decomposition:
\begin{eqnarray}
    \bigcap_{s\geq t\ s\to 
    \infty}\!\!\!\overline{\mathrm{Im}}\bigl[\HHH(\Sigma_{s}, \RR,\mu)\to \HHH(\Sigma_{t}, 
    \RR,\mu)\bigr]=\int^{\oplus}_{X} \bigcap_{s\geq t\ s\to 
    \infty}\!\!\!\overline{\mathrm{Im}}\bigl[\HHH(\Sigma_{s,x})\to 
    \HHH(\Sigma_{t,x})\bigr]d\mu(x)\label{eq: Hilbert integral decomposition}
\end{eqnarray}
    
It remains to make the choice of a $\Sigma$ and of the sequence 
$(\Sigma_{t})_{t\in\Nmath}$ and to relate this with harmonic 
Dirichlet functions via Section~\ref{sect:Harmonic Dirichlet Functions and L2 Cohomology}.
The simplicial complex $\Sigma_{t,x}$ is obtained from $\GG_{x}$ by 
first gluing a disk along each circuit of length  
$t$ and then taking the second barycentric subdivision. The 
simplicial complex $\Sigma_{x}$ is their union.
For each $s$, it follows from naturality in Proposition~\ref{prop:cap H(S-s) to H(S-t) isomorphic with 
HD for simpl-compl}, applied for each $x$, that there is an isomorphism of Hilbert $\RR$-modules:
\begin{eqnarray}
\int^{\oplus}_{X} \bigcap_{s\geq t\ s\to \infty}\!\!\!\overline{\mathrm{Im}}
\bigl[\HHH(\Sigma_{s,x})\to \HHH(\Sigma_{t,x})\bigr]d\mu(x) \simeq \int^{\oplus}_{X} \HD(\GG_{x})/\Cmath \ d\mu(x)\label{eq: 
connection HD-L2}
\end{eqnarray}
so that its $\RR$-dimension
does not depend on $t$, and
Theorem~\ref{thm: Betti 1 and graphings} is proved by putting the 
equalities (\ref{eq:exchanging dim and lim}), (\ref{eq: Hilbert integral decomposition}) and (\ref{eq: 
connection HD-L2}) together.
\endofproof

\section{Some Background about Measured Equivalence Relations}
\label{sect: background Measured Equivalence Relations}

In this section, we just recall briefly the definition of some 
notions appearing in the paper. The reader may consult \cite{FM77}
and \cite{Gab00,Gab02} for more details and more references.

\begin{description}
    \item[Countable standard equivalence relation.] A \emph{countable standard equivalence relation}
    on the standard Borel space $(X,\mu)$ is an equivalence 
    relation $\RR$
    with countable classes that is a Borel subset of $X\times X$ for 
    the product $\sigma$-algebra.
    
    \item[Preservation of the measure.]
    The (countable standard) equivalence relation $\RR$ is said 
    to \emph{preserve the measure}
    if for every partially defined isomorphism 
    $\varphi:{A}\to {B}$ whose graph is contained in 
    $\RR$ ($\{(x,\varphi(x)): x\in A\}\subset \RR$), one has 
    $\mu({A})=\mu({B})$, or equivalently if{f} 
    the measures $\nu_{1}$ and $\nu_{2}$ on the set 
    $\RR\subset X\times X$ coincide, 
    defined with respect to the projections on the first (resp.\  second) coordinate 
    $pr_{1}$ (resp.\  $pr_{2}$) by 
    $\nu_{1}(C)=\int_{X} \#(C\cap pr_{1}^{-1}(x)) d\mu(x)$
    and 
    $\nu_{2}(C)=\int_{X} \#(C\cap pr_{2}^{-1}(y)) d\mu(y)$.
    One denotes by $\nu=\nu_{1}=\nu_{2}$ this common (usually infinite) 
    measure on $\RR$.

    \item[Essentially Free Action.] A Borel action of $\HH$ on a 
    standard probability measure space $(X,\mu)$ is \emph{essentially free}
    if the Borel subset of points $x\in X$ with non-trivial stabilizer 
    ($\mathrm{Stab}_{\HH}(x)=\{h\in \HH: hx=x\}\not=\{id\}$) has 
    $\mu$-measure $0$. The term ``essentially'' is frequently omitted.
    
    \item[Restrictions.]
    Let $(X,\mu)$ be a standard Borel space with a probability 
    measure $\mu$ and $\RR$ a measure-preserving Borel equivalence relation.
    If $Y$ is a Borel subset of $X$ of non-zero measure, denote by 
    $\mu_{Y}:=\frac{\mu_{\vert Y}}{\mu(Y)}$ the normalized probability measure 
    on $Y$. The \emph{restriction} $\RR_{Y}$ of $\RR$ to $Y$ is the 
    $\mu_{Y}$-measure-preserving Borel equivalence relation on $Y$ 
    defined by for every $x,y\in Y$, $x\RR_{Y} y \Leftrightarrow x\RR y$.
    
	\item[Saturation.]
	A Borel subset $U\subset X$ is called $\RR$-\emph{saturated} if it is a
	union of $\RR$-classes. The $\RR$-saturation of a Borel set $U$ is the
	smallest $\RR$-saturated set containing it. It is the union of the 
	$\RR$-classes meeting $U$.	
	
    \item[Finite index.] \label{subsect:background-finite index}
    A sub-equivalence relation ${\mathcal{S}}\subset \RR$ has \emph{finite 
    index} in $\RR$ if each $\RR$-class decomposes into finitely 
    many ${\mathcal{S}}$-classes. If this number is constant, it is 
    called the index of $\mathcal{S}$ in $\RR$ and is denoted by 
    $[\RR:\mathcal{S}]$.
    
    \item[Graphings.] 
    A \emph{ probability measure-preserving oriented graphing} on $(X,\mu)$ 
    is an at most countable family $\Phi=(\varphi_{i})_{i\in I}$ of 
    partial measure-preserving isomorphisms $\varphi_{i}:A_{i}\to 
    B_{i}$ between Borel subsets $A_{i},B_{i}\subset X$.

    A \emph{probability measure-preserving unoriented graphing} $\Psi$ 
    on $(X,\mu)$ is  a Borel 
    subset of $X\times X\setminus\{(x,x):x\in X\}$ that is symmetric under the flip 
    $(x,y)\leftrightarrow (y,x)$ such that the smallest equivalence 
    relation $\RR_{\Psi}$ containing it has coutable classes and is 
    measure-preserving. It provides a Borel choice of pairs 
    of $\RR_{\Psi}$-equivalent points (``neighbors''), and thus a graph structure on each equivalence class of $\RR_{\Psi}$.  When these graphs are (almost) all trees, the graphing is called a \emph{treeing}.
    
    $\RR_{\Psi}$ is \emph{generated} by $\Psi$.
    
    $\Psi$ is a graphing \emph{over} $\RR$ if it is contained in $\RR$.
    
    An oriented graphing defines clearly an unoriented one, by 
    considering the graphs of the $\varphi_{i}, \varphi_{i}^{-1}$'s.
    The terms probability measure-preserving, oriented and  unoriented 
    are frequently omited.

    The notion of \textit{unoriented graphing} has been introduced by S.~Adams in
\cite{Ada90} 
and \textit{oriented graphing} by G.~Levitt, together with the notion of
\textit{cost}, in \cite{Lev95}.

\item[Cost.] The \emph{cost} of an unoriented graphing 
$\Psi$ is the number $\mathrm{cost}(\Psi,\mu):=\frac{1}{2} \nu(\Psi)$, where $\nu$ is the witness measure 
on $\RR_{\Psi}$ for $\RR_{\Psi}$ to preserve the measure $\mu$ of $X$.

The cost of an oriented graphing $\Phi=(\varphi_{i})_{i\in I}$,
is the sum of the measures of the domains $\sum_{i\in I} \mu(A_{i})$.

Except in the obvious cases (redundancy in $\Phi$), the two notions 
coincide. In general the cost of an oriented graphing is greater than 
that of the associated unoriented one.

The {\em cost ($\mathrm{cost}(\RR, \mu)$) of a p.m.p.\  countable
equivalence relation} $\RR$ is the infimum of the costs of the
generating graphings.  The {\em cost of a group} $\GGamma$ 
is the infimum of $\mathrm{cost}(\RR, \mu)$
over all equivalence relations $\RR$ defined by a p.m.p.\ free actions of $\GGamma$
(see \cite{Gab00}).
A comparison has been established between the
cost and the first $L^2$ Betti number: $\beta_{1}(\RR, \mu)\leq \mathrm{cost}(\RR, \mu)-1$
\cite[Cor.~3.23]{Gab02}.  Despite that equality is not known to be
true in general for an $\RR$ with only infinite classes, there is no
(not yet ?) counterexample.
However, when $\Psi$ is a treeing, then $\beta_{1}(\RR_{\Psi}, \mu)=\mathrm{cost}(\RR_{\Psi}, \mu)-1=\mathrm{cost}(\Psi)-1$ (\cite[Cor.~3.23]{Gab02} and \cite[Th. 1]{Gab00}).

\end{description}

\bigskip
\nobreak
\noindent \textsc{D.~G.: UMPA, UMR CNRS 5669, ENS-Lyon, 
69364 Lyon
Cedex 7, FRANCE}

\noindent \texttt{gaboriau@umpa.ens-lyon.fr}

\end{document}